\documentstyle{amsppt} \magnification=\magstep1 \hsize=6 
truein \hcorrection{.375in} \vsize=8.5 truein 
\parindent=20pt \baselineskip=14pt \TagsOnRight
\NoBlackBoxes\footline{\hss\tenrm\folio\hss}

\centerline{The Diagonal Distribution for the Invariant Measure}
\centerline{ of a Unitary Type Symmetric Space}

\vskip.5truein

\centerline{Doug Pickrell}
\centerline{Mathematics Department}
\centerline{University of Arizona}
\centerline{Tucson, Arizona 85721 }
\centerline{pickrell\@math.arizona.edu}

\vskip.5truein 

\flushpar Abstract.  Let $\Theta$ denote an involution for a 
simply connected compact Lie group $U$, let $K$ denote the 
fixed point set, and let $\mu$ denote the $U$-invariant 
probability measure on $U/K$.  Consider the geodesic 
embedding $\phi :U/K\to U:u\to uu^{-\Theta}$ of Cartan.  In this paper 
we compute the Fourier transform of the diagonal 
distribution for $\phi_{*}\mu$, relative to a compatible triangular 
decomposition of $G$, the complexification of $U$.  This 
boils down to a Duistermaat-Heckman exact stationary 
phase calculation, involving a Poisson structure on the 
dual symmetric space $G_0/K$ discovered by Evens and Lu.  

\bigskip

\centerline{\S 0. Introduction.}

\bigskip

Suppose that $K$ is a simply connected compact Lie group, 
and let $G$ denote the complexification.  Given a 
triangular decomposition $\frak g=\frak n^{-}\oplus \frak h\oplus 
\frak n^{+}$, a generic $g\in K$ 
has a unique ``LDU decomposition'', $g=lmau$, where $l\in N^{-}$ 
(lower triangular), $u\in N^{+}$ (upper triangular), and $ma\in H$, 
with $m\in H\cap K$ (unitary) and $a\in exp(\frak h_{\Bbb R}=\frak h
\cap i\frak k)$ 
(positive).  A formula of Harish-Chandra (essentially 
equivalent to the Weyl dimension formula) asserts that 
for $\lambda\in \frak h_{\Bbb R}^{*}$ 
$$\int_Ka^{-i\lambda}=\bold c(2\delta -i\lambda )=\prod_{\alpha >
0}\frac {\langle 2\delta ,\alpha\rangle}{\langle 2\delta -i\lambda 
,\alpha\rangle},\tag 0.1$$
where the integral is with respect to normalized Haar 
measure, the product is over positive complex roots, and 
$2\delta$ is the sum of the positive complex roots (there are 
various ways in which the $\bold c$-function arises, and this 
formula has many extensions and interpretations; 
see e.g. [H2], especially \S 5-6 of chapter IV).  

The purpose of this paper is to present a generalization 
of this formula, and some of the related geometry, in 
which $K$ is replaced by a compact symmetric space.  

Suppose that $X$ is a simply connected compact symmetric 
space with a fixed basepoint.  From this we obtain (1) a 
diagram of groups, 
$$\matrix &&G\\
&\nearrow&&\nwarrow\\
G_0&&&&U\\
&\nwarrow&&\nearrow\\
&&K\endmatrix ,\tag 0.2$$
where $U$ is the universal covering of the identity 
component of $Aut(X)$, $X\simeq U/K$, $G$ is the complexification 
of $U$, and $G_0/K$ is the noncompact type symmetric space 
dual to $X$; and (2) a diagram of equivariant totally 
geodesic (Cartan) embeddings of symmetric spaces:  
$$\matrix U/K&@>{\phi}>>&U\\
\downarrow&&\downarrow\\
G/G_0&@>{\phi}>>&G&{{\psi}\atop {\leftarrow}}&G/U\\
&&\uparrow&&\uparrow\\
&&G_0&{{\psi}\atop {\leftarrow}}&G_0/K\endmatrix .\tag 0.3$$

We also consider one additional ingredient:  a triangular 
decomposition of $\frak g$, $\frak g=\frak n^{-}\oplus \frak h\oplus 
\frak n^{+}$, which is $\Theta$-stable and 
for which $\frak t_0=\frak h\cap \frak k$ is maximal abelian in $
\frak k$, where $\Theta$ is 
the involution corresponding to the pair $(U,K)$.  

Given this triangular structure, a generic element of 
$\phi (U/K)$ can be written as $g=l\bold wa_{\phi}ml^{*\Theta}$, where $
l$, $a_{\phi}$, $m$ 
are roughly as before, and $\bold w\in T_0^{(2)}$ (elements of order 
two); the possible $\bold w$ index the connected components of 
the set of generic elements.  In the special case in 
which $\Theta$ is an inner automorphism, the generalization of 
$(0.1)$ which we consider is of the form 
$$\int_{\phi (U/K)}a_{\phi}^{-i\lambda}=\frac {\vert W(K)\vert}{\vert 
W(U)\vert}\sum_{\bold w}\prod^{\bold w}\frac {\langle\delta ,\alpha
\rangle}{\langle\delta -i\lambda ,\alpha\rangle}\tag 0.4$$
where, given $\bold w$, the product is over positive roots $\alpha$ 
which are of noncompact type for the involution 
$Ad(\bold w)\Theta$, and $\vert W(\cdot )\vert$ denotes the order of the Weyl group.  

The plan of the paper is the following.  In \S 1 we 
compute the intersections of the $\phi$-images in $(0.3)$ with 
the triangular decomposition of $G$.  A notable qualitative 
fact is that just as the map $U/K\to G/G_0$ in $(0.3)$ is a 
homotopy equivalence, so also are the intersections with 
the triangular decomposition of $G$.  Possibly everything 
in this section is known; it can certainly be generalized 
and packaged in various ways (the canonical source is 
[Wolf]).  

In \S 2 the general formulation and a proof of $(0.4)$ is 
presented.  It turns out that, in the inner case for 
example, the $\bold w=1$ term in $(0.4)$ equals 
$$\int_{G_0/K}\bold a(g_0)^{-2\delta -2(\delta -i\lambda )}dV(g_0
K)=\int_{\psi (G_0/K)}a_{\psi}^{-\delta -(\delta -i\lambda )}\tag 0.5$$
where $g_0=\bold l\bold a\bold u$ is an Iwasawa decomposition in $
G$, 
$\psi (g_0K)=g_0g_0^{*}=\bold la_{\psi}\bold l^{*}$ is an LDU decomposition, and the 
integrals are with respect to a $G_0$-invariant measure.  It 
is remarkable, although not a surprise, that 
$\bold a(g_0)^{-2\delta}dV(g_0K)$ is the volume element for a symplectic 
form having a momentum map $log(\bold a(g_0K))$.  Hence $(0.5)$ 
can be evaluated using (a noncompact version of) the 
Duistermaat-Heckman exact stationary phase method.  
The symplectic structure was discovered by Evens and 
Lu, in a general setting ([EL]); the relevance of this 
structure was pointed out to me by Foth and Otto ([FO]), 
to whom I am grateful.  

It is natural to consider the more general integral 
$$\Psi (g)=\int_{G_0}\bold a(g_0g)^{-2\delta -2(\delta -i\lambda 
)}dg_0\tag 0.6$$
for $g\in G_0\backslash G/U$.  This is an eigenfunction for 
$G$-invariant differential operators on $G/U$.  This can 
also be evaluated exactly, by the same method.  
 
In [Pi1,2] I have discussed conjectural generalizations of 
$(0.1)$ and $(0.4)$ to loop spaces, and other kinds of infinite 
symmetric spaces.  The localization argument applies in 
a heuristic way.  In appendix A there is a proof of $(0.1)$, 
involving an explicit factorization of the integral, which 
has elements that seem useful in the loop space context.  

\bigskip

\flushpar Notation.  $\langle\cdot ,\cdot\rangle$ will denote the Killing form for 
$\frak g$.  For an automorphism $\theta$ of $\frak g$, we will often write 
$\theta (x)=x^{\theta}$, and more briefly, $Ad(g)(x)=x^g$.  We will write 
$x=x_{-}+x_{\frak h}+x_{+}$ for the triangular decomposition of $
x\in \frak g$, 
and $x=x_{\frak k}+x_{\frak p}$ for the Cartan decomposition of $
x\in \frak g_0$.  

\bigskip

\centerline{\S 1. Symmetric Spaces and Triangular Decomposition.}

\bigskip

Throughout the remainder of this paper, $U$ will denote a 
simply connected compact Lie group, $\Theta$ will denote an 
involution of $U$, with fixed point set $K$, and $X$ will 
denote the quotient, $U/K$.  This implies that $K$ is 
connected and $X$ is simply connected (Theorem 8.2 of 
[H1]).  

Corresponding to the diagram of groups in $(0.1)$, there is 
a Lie algebra diagram 
$$\matrix &&\frak g=\frak u\oplus i\frak u&&\\
&\nearrow&&\nwarrow&&\\
\frak g_0=\frak k\oplus \frak p&&&&\frak u=\frak k\oplus i\frak p\\
&\nwarrow&&\nearrow\\
&&\frak k\endmatrix \tag 1.1$$
where $\Theta$, acting on the Lie algebra level and extended 
complex linearly to $\frak g$, is $+1$ on $\frak k$ and $-1$ on $
\frak p$.  We let 
$(\cdot )^{-*}$ denote the Cartan involution for the pair $(G,U)$.  
The Cartan involution for the pair $(G,G_0)$ is given by 
$\sigma (g)=g^{-*\Theta}$.  Since $*$, $\Theta$, $\sigma$, and $(
\cdot )^{-1}$ commute, our 
practice of writing $g^{\Theta}$ for $\Theta (g)$, etc, should not cause 
any confusion.  

We have natural maps 
$$\matrix K&\to&U&\to&U/K\\
\downarrow&&\downarrow&&\downarrow\\
G_0&\to&G&\to&G/G_0\endmatrix .\tag 1.2$$
The vertical arrows (given by inclusion) are homotopy 
equivalences; more precisely, there are diffeomorphisms 
(polar or Cartan decompositions) 
$$K\times \frak p\to G_0,\quad U\times i\frak u\to G,\quad U\times_
Ki\frak k\to G/G_0,\tag 1.3$$
in each case given by the formula $(g,X)\to gexp(X)$ (mod 
$G_0$ in the last case).  In turn there are totally geodesic 
embeddings of symmetric spaces 
$$\matrix U/K&@>{\phi}>>&U&:&gK&\to&gg^{-\Theta}\\
\downarrow&&\downarrow\\
G/G_0&@>{\phi}>>&G&:&gG_0&\to&gg^{*\Theta}=gg^{-\sigma}\endmatrix 
,\tag 1.4$$
where the symmetric space structures are derived from 
the Killing form.  

A group element of the form $g=g_1g_1^{-\sigma}$ satisfies the 
equation $g^{*}=g^{\Theta}$ (i.e.  $gg^{\sigma}=1$); $g^{*}=g^{\Theta}$ implies that 
$Ad(g)\circ\sigma$ is an antilinear involution; and if $g=g_1g_1^{
-\sigma}$, 
then $Ad(g)\circ\sigma =Ad(g_1)\circ\sigma\circ Ad(g_1^{-1})$, hence $
\sigma$ and $Ad(g)\circ\sigma$ 
are inner conjugate.  These considerations lead to the 
following well-known 

\proclaim{ (1.5) Proposition} (a) In terms of $g\in G$, 
$$\matrix \phi (U/K)=\{g^{-1}=g^{*}=g^{\Theta}\}_0&\to&U=\{g^{-1}
=g^{*}\}\\
\downarrow&&\downarrow\\
\phi (G/G_0)=\{g^{*}=g^{\Theta}\}_0&\to&G\endmatrix $$
where $\{\cdot \}_0$ denotes the connected component containing 
the identity.  

(b) The connected components of $\{g^{-1}=g^{*}=g^{\Theta}\}$ are 
determined by the map which sends $g$ to the inner 
conjugacy class of the involution $\eta =Ad(g)\circ\Theta$, subject 
to the constraint that $\eta$ equals $\Theta$ in 
$Out(U)=Ad(U)\backslash Aut(U)$.  A similar statement applies to 
$\{g^{*}=g^{\Theta}\}$, with $\sigma$ and antilinear automorphisms of $
G$ in 
place of $\Theta$ and involutions of $K$.  
\endproclaim

\demo{Proof of (1.5)} We first recall why $\{gg^{\sigma}=1\}$ is 
smooth.  

Consider the map $\psi :G\to G:g\to gg^{\sigma}$.  If we use right 
translation to identify the tangent space at any point of 
$G$ with $\frak g$, the derivative at $g$ is given by 
$x\to x+Ad(g)[\sigma (x)]$.  Thus $ker(d\psi\vert_g)$ is identified with the 
$-1$ eigenspace of $Ad(g)\circ\sigma$ acting on $\frak g$.  

Now suppose $gg^{\sigma}=1$.  Since $Ad(g)\circ\sigma$ is an involution, 
the spectrum of $Ad(g)\circ\sigma$ is fixed.  Thus the dimension 
of the $-1$ eigenspace of $Ad(g)\circ\sigma$ is constant on $\{gg^{
\sigma}=1\}$.  
It follows that $\psi$ has constant rank on the connected 
components of $\psi^{-1}(1)$.  Since $\psi$ is an algebraic map, this 
implies that $\{g^{*}=g^{\Theta}\}$ is an embedded submanifold.  A 
similar argument applies to the intersection with $U$.  

The action 
$$G\times \{gg^{\sigma}=1\}\to \{gg^{\sigma}=1\}:g,g_1\to gg_1g^{
*\Theta}\tag 1.6$$
is isometric (for the symmetric space structure).  The 
constancy of the rank of $\psi$ on connected components is 
equivalent to the statement that the dimension of the 
isotropy subgroup for the action of $G$ is constant on 
connected components of $\{gg^{\sigma}=1\}$ (in fact this dimension 
is the same on all components).  Hence the action of $G$ 
must be transitive on connected components.  The same 
applies to the same action of $U$ on $\{g\in U:gg^{\Theta}=1\}$.  This 
implies $(a)$.  

For the first part of $(b)$, note that in fact the map 
$$\{g\in U:g^{-1}=g^{\Theta}\}\to \{\eta\in Aut(U)^{(2)}:Out(\eta 
)=Out(\Theta )\}:g\to Ad(g)\circ\Theta$$
is a universal covering for each connected component 
(For the identity component this covering is understood 
more intellibly by identifying the total space with $U/K$:  
$$C_U(K)/C(U)\to U/K@>{q}>>Ad(U)\cdot\Theta\tag 1.7$$
where $q(g_1K)=Ad(g_1)\circ\Theta\circ Ad(g_1)^{-1}$; to obtain a similar 
picture for another component, we replace $\Theta$ by 
$Ad(g)\circ\Theta$, for some $g$ in the component).  

The second part of $(b)$ is similar (We could also note 
that the inclusion $\{g^{-1}=g^{*}=g^{\Theta}\}\to \{g^{*}=g^{\Theta}
\}$ is a 
homotopy equivalence, since we know this is true for 
the identity component, and we are free to change $\Theta$ to 
$Ad(g)\circ\Theta$; the fact that the $\pi_0$'s are the same is a 
reflection of the fact that classifying $\Theta$'s and classifying 
$\sigma$'s are canonically isomorphic problems (see e.g.  2.  of 
\S 6, chapter 10 of [H])).   \qed
\enddemo\ 

\flushpar Remarks (1.8).  (a) I do not know of a uniform 
way to define an invariant for the class of an involution 
$\eta$ as in $(b)$.  However it is a simple matter to produce 
an invariant in a case by case manner from the 
classification of symmetric spaces (see Table V of [H1]).  

(b) The groups and maps in $(1.4)$ exist for any 
automorphism $\Theta$ of $K$.  However it seems that there is a 
linear characterization of the $\phi$-images (up to 
connectedness issues), and $G_0$ is a real form, only in the 
symmetric case, $\Theta^2=1$.  

\smallskip

Fix a maximal abelian subalgebra $\frak t_0\subset \frak k$.  We then obtain 
$\Theta$-stable Cartan subalgebras 
$$\frak h_0=\Cal Z_{\frak g_0}(\frak t_0)=\frak t_0\oplus \frak a_
0,\quad \frak t=\frak t_0\oplus i\frak a_0,\quad and\quad \frak h
=\frak h_0^{\Bbb C}\tag 1.9$$
for $\frak g_0$, $\frak u$, and $\frak g$, respectively, where $\frak a_
0\subset \frak p$ (see $(6.60)$ of 
[Kn]).  We let $T_0$ and $T$ denote the maximal tori in $K$ 
and $U$ corresponding to $\frak t_0$ and $\frak t$, respectively.  

Let $\Delta$ denote the roots for $\frak h$ acting on $\frak g$; $
\Delta\subset \frak h_{\Bbb R}^{*}$, where 
$\frak h_{\Bbb R}=\frak a_0\oplus i\frak t_0$.  We choose a Weyl chamber $
C^{+}$ which is 
$\Theta$-stable (to prove that $C^{+}$ exists, we must show that 
$i\frak t_0$, the $+1$ eigenspace of $\Theta$ acting on $\frak h_{
\Bbb R}$, contains a 
regular element of $\frak g$; this is equivalent to the fact that 
$\frak h_0$ in $(1.8)$ is a Cartan subalgebra).  Since $\sigma =-
(\cdot )^{*\Theta}$ and 
$(\cdot )^{*}$ is the identity on $\frak h_{\Bbb R}$, $\sigma (C^{
+})=-C^{+}$.  

Given our choice of $C^{+}$, we obtain a $\Theta$-stable triangular 
decomposition $\frak g=\frak n^{-}\oplus \frak h\oplus \frak n^{+}$, so that $
\sigma (\frak n^{\pm})=\frak n^{\mp}$.  Let 
$N^{\pm}=exp(\frak n^{\pm})$, $H=exp(\frak h)$, and $B^{\pm}=HN^{
\pm}$.  We also let 
$W=W(G,T)$ denote the Weyl group, 
$W=N_U(T)/T\simeq N_G(H)/H$.  

At the group level we have the Birkhoff or triangular or 
LDU decomposition for $G$, 
$$G=\bigsqcup_W\Sigma^G_w,\quad\Sigma^G_w=N^{-}wHN^{+},\tag 1.10$$
where $\Sigma_w^G$ is diffeomorphic to $(N^{-}\cap wN^{-}w^{-1})\times 
H\times N^{+}$.  
When we intersect this decomposition with $\phi (G/G_0)$, and 
the other spaces in (b) of $(1.5)$, we obtain various 
decompositions.  We will first determine the structure 
of the pieces in the $\{g^{*}=g^{\Theta}\}$ case (thus we initially
ignore connectedness issues).  

\proclaim{ (1.11) Proposition} Fix $w\in W$.  

(a) The intersection $\{g^{*}=g^{\Theta}\}\cap\Sigma^G_w$ is nonempty if and 
only if there exists $\bold w\in w\subset N_U(T)$, such that $\bold w^{
*\Theta}=\bold w$; 
$\bold w$ is unique modulo the action 
$$T\times \{\bold w\in N_U(T):\bold w^{*\Theta}=\bold w\}\to \{\bold w^{
*\Theta}=\bold w\}:\lambda ,\bold w\to\lambda \bold w\lambda^{*\Theta}
.$$

(b) For the action $B^{-}\times \{g^{*}=g^{\Theta}\}\to \{g^{*}=g^{
\Theta}\}$:$b,g\to bgb^{*\Theta}$, 
the stability subgroup is given by 
$$B^{-}_{\bold w}=\{b:\bold w^{-1}b\bold w=\sigma (b)\}\simeq \{l
\in N^{-}:\bold w^{-1}l\bold w=\sigma (l)\in N^{+}\}\times \{h\in 
H:h^{w^{-1}}=\sigma (h)\}.$$

(c) The orbits of $B^{-}$ in $\{g^{*}=g^{\Theta}\}\cap\Sigma^G_w$ are open and 
indexed by 
$$\pi_0(\{\bold w\in w:\bold w^{*\Theta}=\bold w\})\simeq \{\bold w
\in w:\bold w^{*\Theta}=\bold w\}/T,$$
where $T$ acts as in part (a).  

(d) The map 
$$N^{-}\cap N^{-w}\times \{h\in H,L\in N^{-}\cap N^{+w}:h^{w^{-1}
\Theta}=h^{*},\sigma (L)^{\bold wh}=L^{-1}\}\to \{g^{*}=g^{\Theta}
\}\cap\Sigma^G_w$$
given by $l,h,L\to lL^{-1}\bold wh(lL^{-1})^{*\Theta}$ is a diffeomorphism onto 
the connected component containing $\bold w$.  This component 
is homotopic to the torus $exp(\{Ad(w^{-1})\Theta\vert_{\frak t}=
-1\})$.  

(e) In particular for $w=1$, the map 
$$N^{-}\times (T_0^{(2)}\times_{exp(i\frak a_0)^{(2)}}exp(i\frak a_
0))\times exp(i\frak t_0)\to \{g^{*}=g^{\Theta}\}\cap\Sigma^G_1$$
$$l,[\bold w,m],a_{\phi}\to g=l\bold wma_{\phi}l^{*\Theta}$$
is a diffeomorphism, so that the connected components 
for $\{g^{*}=g^{\Theta}\}\cap\Sigma_1^G$ are indexed by $T_0^{(2)}
/exp(i\frak a_0)^{(2)}$.  

\endproclaim

\demo{Proof of (1.11)} Suppose that $g\in\Sigma^G_w$.  We write 
$g=l\bold whu$, for some $l\in N^{-}$, $\bold w\in w\subset N_U(T
)${\bf ,} $h\in exp(\frak h_{\Bbb R})$, 
$u\in$$N^{+}$.  If we additionally require that $l\in N^{-}\cap (
N^{-})^w$, 
then this decomposition is unique, but we will not 
require this at the outset.  

We have $g=g^{*\Theta}$ if and only if 
$$l\bold whu=u^{*\Theta}(\bold wh)^{*\Theta}l^{*\Theta}\tag 1.12$$
if and only if 
$$(\bold wh)^{*\Theta}=(u^{\sigma}l)(\bold wh)(ul^{\sigma})$$
$$=\{(u^{\sigma}l)_{-}(u^{\sigma}l)_{+}\}(\bold wh)(ul^{\sigma})$$
$$=(u^{\sigma}l)_{-}(\bold wh)\{(u^{\sigma}l)_{+}^{(\bold wh)^{-1}}
ul^{\sigma}\}\tag 1.13$$
where $L=L_{-}L_{+}$ denotes the decomposition induced by the 
diffeomorphism 
$$N^{-}\cap_{}(N^{-})^w\times_{}N^{-}\cap (N^{+})^w\to N^{-}:L_{-}
,L_{+}\to L=L_{-}L_{+}\tag 1.14$$
Thus $(1.13)$ holds if and only if 
$$(u^{\sigma}l)_{-}=1=(u^{\sigma}l)_{+}^{(\bold wh)^{-1}}ul^{\sigma}\tag 1.15$$
and $(\bold wh)^{*\Theta}=\bold wh$, or, using the fact that $h$ is real, 
$$h^{\Theta}=h^w\quad and\quad \bold w\bold w^{\Theta}=1.\tag 1.16$$

Consider part (a).  If $g$ is in the intersection, then we 
have just seen that $\bold w$ must satisfy $\bold w^{*\Theta}=\bold w$.  
Conversely, given a unitary representative $\bold w$ for $w$ 
satisfying $\bold w^{*\Theta}=\bold w$, the intersection contains $
\bold w$ and 
hence is nonempty.  This proves (a).  

Part (b) is straightforward.  

Now consider (c).  We first write $g\in \{g^{*}=g^{\Theta}\}\cap\Sigma^
G_w$ 
uniquely as $l\omega u$, where $l\in N^{-}\cap wN^{-}w^{-1}$, $\omega 
=\bold wh$, and 
$u\in N^{+}$.  For the first part of (c) we must prove that 
we can relax the constraint on $l$ to arrange for $u=l^{*\Theta}$.  
We can write 
$$g=lL^{-1}\omega \{(\omega L\omega^{-1})u\},\tag 1.17$$
where $L\in N^{-}\cap wN^{+}w^{-1}$ is arbitrary.  We must prove the 
existence of $L$ such that $lL^{-1}=u^{*\Theta}\omega\sigma (L^{-
1})\omega^{-1}$, or 
$$u^{\sigma}l=\omega L^{-\sigma}\omega^{-1}L.\tag 1.18$$

The basic fact is that this equation has a unique 
solution $L\in N^{-}\cap N^{+w}$ satisfying $\omega L^{\sigma}\omega^{
-1}=L^{-1}$, namely 
$L=(u^{\sigma}l)^{1/2}$ (square root has an unambiguous meaning in 
a simply connected nilpotent Lie group).  To see this 
simply plug such an $L$ into $(1.18)$.  We obtain the 
equation $L^2=u^{\sigma}l$.  The fact that $u^{\sigma}l$, and its square 
root, satisfy $\omega L^{\sigma}\omega^{-1}=L^{-1}$ follows from $
(1.15)$, and 
uniqueness of the square root.  

As we remarked previously, the existence of a solution $L$ 
proves that $B^{-}$ has open orbits.  The rest of part (c) is 
relatively straightforward, using $(1.16)$.  

The uniqueness of the solution $L$, subject to the 
constraint we imposed, implies the first part of (d).  
The second statement in $(d)$ follows routinely from the 
first part.  

For part (e), to clarify the statement, observe that 
$$exp(i\frak a_0)^{(2)}=K\cap exp(i\frak a_0)=T_0^{(2)}\cap exp(i
\frak a_0).\tag 1.19$$
Now suppose that $w=1$.  In this case $\bold w$ is in the kernel 
of the homomorphism $T\to T:\bold w\to \bold w\bold w^{\Theta}$, and this equals the 
subgroup generated by $T_0^{(2)}$ and $exp(i\frak a_0)$.  We can modify 
$\bold w$ by multiplying by something in the image of the 
homomorphism $T\to T:\lambda\to\lambda\lambda^{-\Theta}$.  This image is $
exp(i\frak a_0)$.  
Therefore we can choose $\bold w\in T_0^{(2)}$, but this choice is 
unique only modulo the intersection of $T_0^{(2)}$ and $exp(i\frak a_
0)$.  
This proves (e).  \qed
\enddemo

\flushpar Example (1.20).  In the group case, $U=K\times K$, 
where $K$ embeds diagonally.  The image of $\frak t_0$ inside $\frak u$ is 
$\{(x,x):x\in \frak t_0\}$, while $i\frak a_0=\{(x,-x):x\in \frak t_
0\}$.  This implies that 
the quotient $T_0^{(2)}/exp(i\frak a_0)^{(2)}$ is trivial.  Thus in this 
group case, the set of generic elements (considered in 
part (e)) is connected, as we already know.  

\proclaim{ (1.22) Notation}Given $\bold w$ as in $(c)$ of $(1.11)$, we 
let $\Sigma^{\{g^{*}=g^{\Theta}\}}_{\bold w}$ denote the corresponding connected 
component of $\{g^{*}=g^{\Theta}\}\cap\Sigma^G_w$ (the $B^{-}$-orbit of $
\bold w$, in the 
sense of (c) of $(1.11)$).  If $\bold w\in\phi (G/G_0)$, then we will 
write $\Sigma_{\bold w}^{\phi (G/G_0)}$ for this component.  We also set 
$\Sigma_{\bold w}^{\phi (U/K)}=\phi (U/K)\cap\Sigma_{\bold w}^{\phi 
(G/G_0)}$.  
\endproclaim

Having understood the intersection of $\{g^{*}=g^{\Theta}\}$ with the 
triangular decomposition, we now want to specialize this 
to the identity component.  In an abstract way this is 
answered by $(1.5b)$ (and Remark $(1.8a)$).  Concerning open 
orbits, we have the following 

\proclaim{ (1.23) Proposition} Suppose that $\bold w\in N_U(T)$ 
satisfies $\bold w^{-\Theta}=\bold w$.  The following are equivalent:  

(a) $\Sigma_{\bold w}^{\{g^{*}=g^{\Theta}\}}$ is an open $B^{-}$-orbit in the identity 
component, $\phi (G/G_0)$.  

(b) There exists $\bold w_1\in N_U(T_0)$ such that $\phi (\bold w_
1K)=\bold w$.  

Hence the open orbits can be parameterized by either 
$N_U(T_0)/N_K(T_0)$ (the intrinsic point of view), or the set 
of $\bold w\in T_0^{(2)}/exp(i\frak a_0)^{(2)}$ such that $Ad(\bold w
)\circ\Theta$ is equivalent 
to $\Theta$ in the sense of $(1.5b)$ (the nonintrinsic point of 
view, as in $(1.11e)$).  

In addition, the $\bold w_1K$ are exactly the $T_0$ fixed points in 
$U/K$.  
\endproclaim

\demo{Proof of (1.23)}Determining the possible (open) $B^{-}$ 
orbits in $G/G_0$ is equivalent to determining the possible 
(open) $G_0$ orbits in $B^{-}\backslash G$.  Thus the equivalence of $
(a)$ 
and $(b)$ follows from 
Theorem 4.6 and its Corollaries in [Wolf]. The other
statements are obvious.
\qed
\enddemo
 
In general it apparently remains an open question to 
systematically obtain representatives for all $B^{-}$ orbits in 
$G/G_0$, from the intrinsic point of view (see [WZ] for 
the Hermitian symmetric case).  In this regard the 
nonintrinsic point of view of $(1.11)$ seems to have some 
utility.  

Let $q:G\to G/B^{+}$ denote the quotient map.  The map $q$ 
applied to the decomposition $(1.9)$ induces the (more 
conventional) triangular stratification for the flag space, 
$$U/T\simeq G/B^{+}=\bigsqcup_W\Sigma_w,\quad\Sigma_w=N^{-}\cdot 
wB^{+},\tag 1.24$$
where each $\Sigma_w$ is a cell ($\simeq N^{-}\cap wN^{-}w^{-1}$).  As a 
consequence, for the pieces of the induced decomposition 
for $U$, there are diffeomorphisms 
$$\Sigma_w^U=U\cap\Sigma^G_w\simeq\Sigma_w\times T.\tag 1.25$$
The inclusions $\Sigma^U_w\to\Sigma^G_w$ are homotopy equivalences, 
because $T$ is homotopy equivalent to $B^{+}$:  
$$\matrix T&\to&\Sigma_w^U&@>{q}>>&\Sigma_w\\
\downarrow&&\downarrow&&\Vert\\
B^{+}&\to&\Sigma_w^G&@>{q}>>&\Sigma_w\endmatrix \tag 1.26$$

The main point of this section is now to describe the 
generalization of this to $U/K\to G/G_0$.  We consider the 
Iwasawa decomposition for $G$, which we write as 
$$G\simeq N^{-}\times A\times U:g=\bold l(g)\bold a(g)\bold u(g),\tag 1.27$$
where $A=exp(\frak h_{\Bbb R})$.  There is an induced right action 
$$U\times T\times G_0\to U:(u,t,g_0)\to t^{-1}\bold u(ug_0)\tag 1.28$$
arising from the identification of $U$ with $N^{-}A\backslash G$.  

\proclaim{ (1.29) Proposition} Suppose that $\bold w\in N_U(T)$ 
satisfies $\bold w^{-\Theta}=\bold w$ and $\Sigma^{\{g^{*}=g^{\Theta}
\}}_{\bold w}\subset\phi (G/G_0)$.  Fix a choice 
of $\bold w_1\in U$ such that $\bold w_1\bold w_1^{-\Theta}=\bold w$.  

(a) The map 
$$T\times G_0\to\Sigma^{\phi (U/K)}_{\bold w}:(t,g_0)\to\phi (t^{
-1}\bold u(\bold w_1g_0))\tag 1.30$$
is surjective, and induces a diffeomorphism 
$$T\times_{exp(\{Ad(\bold w)\Theta\vert_{\frak t}=1\})}R\backslash 
G_0/K\to\Sigma^{\phi (U/K)}_{\bold w}\tag 1.31$$
where $R=(N^{-}A)^{\bold w_1^{-1}}\cap G_0$ is a contractible subgroup of 
$G_0$, and $\lambda\in exp(\{Ad(\bold w)\Theta\vert_{\frak t}=1\}
)$ is identified with a pair 
$(\lambda ,\lambda^{\bold w_1^{-1}})$.  

(b) The inclusion 
$$\Sigma^{\phi (U/K)}_{\bold w}\to\Sigma^{\phi (G/G_0)}_{\bold w}
.\tag 1.32$$
is a homotopy equivalence; each is homotopic to the 
torus $exp(\{Ad(w^{-1})\Theta\vert_{\frak t}=-1\})$.  

(c) The connected components of $\phi (U/K)$ intersected with 
$\Sigma^U_1$ are indexed by $\bold w\in T_0^{(2)}/exp(i\frak a_0)^{
(2)}$ such that 
$\bold w=\bold w_1\bold w_1^{-\Theta}$, for some $\bold w_1\in N_
U(T_0)$; for such a $\bold w$, and 
choice of $\bold w_1$, the diffeomorphism in (a) simplifies to 
$$exp(i\frak a_0)/exp(i\frak a_0)^{(2)}\times A_0\backslash G_0/K
\to\Sigma_{\bold w}^{\phi (U/K)},\tag 1.33$$
where $A_0=exp(\frak a_0)$.  
\endproclaim

\demo{Proof of $(1.29)$} In proving the first part of (a), it 
is convenient to work with $U/K$ instead of $\phi (U/K)$.  
Thus we consider a point in the intersection of $U/K$ and 
the $B^{-}$-orbit of $\bold w_1G_0\in G/G_0$.  This point can be 
represented by a $u\in U$ such that $u=b^{-}\bold w_1g_0$, for some 
$b^{-}\in B^{-}$ and $g_0\in G_0$.  This immediately implies that $
u$ is 
in the $T\times G_0$-orbit of $\bold w_1$, and this proves surjectivity of 
the first map in (a).  

For the second part of (a), we first calculate the 
stabilizer for the action $(1.28)$ at the point $\bold w_1$.  Suppose 
that $t\in T$ and $g_0\in G_0$ satisfy $t\bold u(\bold w_1g_0)=\bold w_
1$.  This is 
equivalent to 
$$\bold w_1g_0\bold w_1^{-1}=\bold l(\bold w_1g_0)\bold a(\bold w_
1g_0)t.\tag 1.34$$
This implies that $g_0\in G_0\cap (B^{-})^{\bold w_1^{-1}}$, and $
t=T(g_0^{\bold w_1})$.  
Conversely if $g_0\in G_0\cap (B^{-})^{\bold w_1^{-1}}$, then $(1
.34)$ holds with 
$t=T(g_0^{\bold w_1})$.  Thus the stabilizer is isomorphic to 
$$\{(T(g_0^{\bold w_1}),g_0):g_0\in G_0\cap (B^{-})^{\bold w_1^{-
1}}\}\subset T\times G_0.\tag 1.35$$

The group $G_0\cap (B^{-})^{\bold w_1^{-1}}$ is connected and solvable.  The 
torus part is isomorphic to $\{\lambda\in T:\lambda^{\bold w_1^{-
1}}\in G_0\}$.  This 
condition on $\lambda$ is equivalent to $(\lambda^{\bold w_1^{-1}}
)^{\Theta}=\lambda^{\bold w_1^{-1}}$, or 
$\lambda\in exp(\{Ad(\bold w)\Theta\vert_{\frak t}=1\})$.  This implies $
(1.31)$.  

From (a), since $R\backslash G_0/K$ is contractible, it follows that 
the double coset space in $(1.31)$ is homotopic to 
$exp(\{Ad(\bold w)\Theta\vert_{\frak t}=-1\})$, modulo elements of order $
2$.  A $t$ in 
this torus is mapped in $(1.30)$ to $t^{-1}\bold wt^{\Theta}=\bold w
(t^{-\bold w^{-1}}t^{\Theta})$.  It 
is straightforward to check that $t^{-\bold w^{-1}}t^{\Theta}$ belongs to 
$exp(\{Ad(\roman w^{-1})\Theta\vert_{\frak t}=-1\})$.  Together with $
(d)$ of $(1.11)$, this 
implies $(b)$.  

Part $(c)$ follows from $(a)$.  \qed
\enddemo

We want to explain how this proposition is related to 
familiar facts in special cases.  First consider the group 
case $X=K$.  We have already explained why the generic 
set is connected; see Example $(1.20)$.  In the group case, 
$(1.1)$ has the form 
$$\matrix &&\{(x,y)\in \frak k^{\Bbb C}\oplus \frak k^{\Bbb C}\}\\
&\nearrow&&\nwarrow\\
\{(x,-x^{*}):x\in \frak k^{\Bbb C}\}&&&&\{(x,y)\in \frak k\oplus 
\frak k\}\\
&\nwarrow&&\nearrow\\
&&\{(x,x):x\in \frak k\}\endmatrix \tag 1.36$$
Given $g_0\in G_0=K^{\Bbb C}$, $g_0$ maps to $(g_0,g_0^{-*})\in G
=K^{\Bbb C}\times K^{\Bbb C}$.  
Given an arbitrary triangular decomposition 
$\frak k^{\Bbb C}=\tilde {\frak n}_{-}\oplus\tilde {\frak h}\oplus
\tilde {\frak n}_{+}$, we obtain a $\Theta$-stable triangular 
decomposition for $\frak g$ by defining $\frak n_{\pm}=\tilde {\frak n}_{
\pm}\times\tilde {\frak n}_{\pm}$, and 
$\frak h=\tilde {\frak h}\times\tilde {\frak h}$.  The Iwasawa factorization $
(1.27)$ is equivalent 
to the two Iwasawa factorizations 
$$g_0=l_1a_1k_1,\quad g_0^{-*}=l_2a_2k_2,\tag 1.37$$
where $a_i\in exp(\tilde {\frak h}_{\Bbb R})$, $l_i\in\tilde {N}^{
-}$, $k_i\in K$.  We have an 
equivariant isomorphism $U/K\to K:(g,h)\to gh^{-1}$.  The map in 
(c) of $(1.29)$ (using $(\cdot )^{*}$ in place of inverse) is given by 
$$T/T^{(2)}\times A\backslash K^{\Bbb C}/K\to\Sigma^K_1:t,g_0K\to 
ta_1^{-1}l^{-1}_1l_2^{-*}a_2^{-*}t.\tag 1.38$$
So $a_{\phi}=a_1^{-1}a_2^{-*}=(a_1a_2)^{-1}$.  Thus in $(1.29)$ we are using 
$exp(\tilde {\frak h}_{\Bbb R})g_0K\in exp(\tilde {\frak h}_{\Bbb R}
)\backslash K^{\Bbb C}/K$ as coordinate, which is 
completely equivalent to using $l_1$ or $l_2\in N^{-}$.  From this 
point of view, $l_1$ is a horocycle coordinate.  

Now suppose that $U/K$ is Hermitian symmetric.  In this 
case $\Theta$ is inner, $P=K^{\Bbb C}B^{+}$ is a parabolic subgroup of $
G$, 
and the natural map $\frak i:U/K\to G/P$ is a $U$-equivariant 
isomorphism.  The natural map $\eta :G_0/K\to G/P$ is an open 
holomorphic embedding, and the image is contained 
$\frak i((U/K)_r)$, the regular set (see ch VIII of [H1], especially 
Prop 7.14).  

There is a commutative diagram 
$$\matrix G_0/K&@>{\bold u}>>&\Sigma^{U/K}_1:&g_0K\to \bold u(g_0
)K&\\
\downarrow I&&\downarrow \frak i\\
G_0/K&@>{\eta}>>&G/P\endmatrix \tag 1.39$$
where the top arrow $\bold u$ is a diffeomorphism, by $(b)$ of 
$(1.29)$, and the map $I$ is defined in the following way:  
given $g_0K\in G_0/K$, we can write $\bold u(g_0)=exp(ix)k$ 
uniquely, where $x\in \frak p$, $exp(itx)K$ is a geodesic of minimal 
length joining the basepoint to $\bold u(g_0)K$, and $k\in K$; we set 
$I(g_0K)=(g_0k^{-1})^{-1}K$.  To see that the diagram is 
commutative, note that because $g_0k^{-1}\in G_0$, 
$(g_0k^{-1})^{-1}=(g_0k^{-1})^{*\Theta}=exp(ix)a^{*\Theta}l^{*\Theta}$, and $
l^{*\Theta}\in N^{+}$; thus 
$I(g_0K)$ equals $\bold u(g_0)$ mod $P$.  

Thus in the Hermitian symmetric case, $\Sigma_1^{U/K}$ is the 
usual model of $G_0/K$ inside $U/K$, but the 
parameterization in $(b)$ of $(1.29)$ is related to the natural 
holomorphic map $\eta$ in a clumsy way.  

\bigskip

\centerline{\S 2. Diagonal Distribution}

\bigskip

Suppose that $\bold w\in T_0^{(2)}=T\cap \{g^{*}=g^{\Theta}\}$.  Then $
Ad(\bold w)\Theta$ is an 
involution, and $\frak n^{-}\oplus \frak h\oplus \frak n^{+}$ is $
Ad(\bold w)\Theta$-stable.  We also 
suppose that $\bold w\in\phi (U/K)$.  

Given $g\in\Sigma^{\phi (U/K)}_{\bold w}$, we factor $g$ as in (e) of $
(1.11)$, 
$$g=l\bold wma_{\phi}l^{*\Theta},\tag 2.1$$
where $l\in N^{-}$, $a_{\phi}\in exp(i\frak t_0)$, and 
$[\bold w,m]\in T_0^{(2)}\times_{exp(i\frak a_0)^{(2)}}exp(i\frak a_
0)$.  

\proclaim{ (2.2) Theorem} For $\lambda\in (i\frak t_0)^{*}$, 
$$\int_{\Sigma^{\phi (U/K)}_{\bold w}}a_{\phi}(g)^{-i\lambda}=\frac 
1M\prod^{\bold w}\frac {\langle\delta ,\alpha\rangle}{\langle\delta 
-i\lambda ,\alpha\rangle}\tag 2.3$$
where the product is over pairs $(\alpha ,Ad(\bold w)\Theta (\alpha 
))$ of 
positive complex roots which are $\underline {not}$ of compact type 
for $Ad(\bold w)\Theta$, and $M=\vert N_U(T_0)/N_K(T_0)\vert$.  Hence 
$$\int_{\phi (U/K)}a_{\phi}(g)^{-i\lambda}=\frac 1M\sum (\prod^{\bold w}\frac {
\langle\delta ,\alpha\rangle}{\langle\delta -i\lambda ,\alpha\rangle}
)\tag 2.4$$
where the sum is over representatives $\bold w$ for the 
connected components of $\phi (U/K)\cap\Sigma^U_1$.  
\endproclaim

Note that it does not matter whether we take $\alpha$ or 
$Ad(\bold w)\Theta (\alpha )$ in the product $(2.3)$, because $\delta$ and $
\lambda$ are fixed 
by $Ad(\bold w)\Theta$.  In the case in which $\Theta$ is inner, i.e.  
$\frak a_0=0$, all roots are either of compact or noncompact 
type.  Hence in this case the product in $(2.3)$ is over 
the noncompact type roots.  

There is a more intrinsic way to write $(2.4)$.  The right 
hand side can be expressed as a sum over 
$\bold w_1\in N_U(T_0)/N_K(T_0)$, using $Ad(\bold w)\Theta =Ad(\bold w_
1)\Theta Ad(\bold w_1^{-1})$ (see 
$(1.23)$).  

To prove $(2.2)$ we will first 
note that we can reduce to the case $\bold w=1$ in $(2.3)$.  
We will 
then need several Lemmas.  

To see that it suffices to prove $(2.3)$ in the case $\bold w=1$, 
observe that $\underline {left}$ translation by $\bold w\in T^{(2
)}_0$, which is an 
isometric map for the Riemannian structure of $U$, maps 
the $\bold w=1$ component to the $\bold w$-component:  
$$L_{\bold w}:\Sigma^{\phi (U/K)}_1\to\Sigma^{\phi (U/K)}_{\bold w}
.\tag 2.5$$
This map interchanges the canonical factorization from 
that relative to $\Theta$ to the one relative to $Ad(\bold w)\Theta$:  if $
g$ 
has the unique decomposition $g=lma_{\phi}l^{*\Theta}$, then $\bold w
g$ has 
the unique decomposition $\bold wg=l^{\bold w}\bold wma_{\phi}(l^{
\bold w})^{*Ad(\bold w)\Theta}$.  
Since $a_{\phi}$ is unchanged, the integral is evaluated in the 
same way, except the meaning of the roots changes.  

We henceforth suppose $\bold w=1$.  Consider the 
parameterization 
$$\Phi :exp(i\frak a_0)\times A_0\backslash G_0/K\to\Sigma_1^{\phi 
(U/K)}:(t,A_0g_0K)\to\phi (t^{-1}\bold u(g_0)).\tag 2.6$$
from $(1.33)$.  Note that $a_{\phi}=\bold a(g_0)^{-1}\bold a(g_0)^{
\tau}\in exp(i\frak t_0)$.  

\proclaim{ (2.7) Lemma} We have 
$$\Phi^{*}(dV_{U/K})=a_{\phi}^{2\delta}(dV_{exp(i\frak a_0)}\times 
dV_{A_0\backslash G_0/K})$$
where $dV_{A_0\backslash G_0/K}$ is obtained by integrating a 
$G_0$-invariant measure on $A_0\backslash G_0$ over $K$.  Thus $(
2.3)$ 
equals 
$$\int a_{\phi}(\bold u(g_0))^{2\delta -i\lambda}dV_{A_0\backslash 
G_0/K}.\tag 2.8$$
\endproclaim\ 
 
\demo{Proof of (2.7)} We will consider a slight 
reformulation of the problem.  We identify $U/K$ with 
$\phi (U/K)$.  Let $S$ denote the inverse image of $\Sigma_1^{U/K}$ in $
U$, 
with respect to the projection $U\to U/K$.  Consider the 
lift 
$$\Psi :exp(i\frak a_0)\times A_0\backslash G_0\to S:(t,A_0g_0)\to 
t^{-1}\bold u(g_0)\tag 2.9$$
of $\Phi$ in $(2.6)$.  We must show that the Jacobian for the 
mapping $\Psi$, with respect to the Riemannian structures 
induced by the Killing form, is equal to a constant times 
$a_{\phi}^{2\delta}$.  To do this we identify $i\frak a_0$, $\frak g_
0\ominus \frak a_0$, and $\frak u$ with the 
tangent spaces to $exp(i\frak a_0)$, $A_0\backslash G_0$ and $U$, respectively, 
using the exponential map and right translation (we use 
right translation because $A_0$ appears on the left).  Let 
$P:\frak g\to \frak u$ denote the projection with kernel $\frak n^{
-}\oplus \frak h_{\Bbb R}$.  We 
compute 
$$d\Psi\vert_{(t,A_0g_0)}:i\frak a_0\oplus (\frak g_0\ominus \frak a_
0)\to \frak u:(\chi ,x)\to\frac d{d\epsilon}\vert_{\epsilon =0}(t
e^{\epsilon\chi})^{-1}\bold u(e^{\epsilon x}g_0)\bold u(g_0)^{-1}
t$$
$$=Ad(t^{-1})\{-\chi +P(Ad(\bold a^{-1}\bold l^{-1})(x))\}.\tag 2.10$$

The operator $Ad(t^{-1})$ preserves $\frak u$-volume, so it can be 
ignored.  

Write $\bold a=\bold a_1\bold a_0$, relative to the decomposition 
$A=exp(i\frak t_0)A_0$.  Since $\bold a_0\in G_0$, $Ad(\bold a_0)$ will preserve 
$\frak g_0$-volume.  Thus the determinant of $(2.10)$ equals the 
determinant of the map 
$$\frak g_0\ominus \frak a_0\to \frak u\ominus i\frak a_0:x\to P_
1(Ad(\bold l')Ad(\bold a_1^{-1})(x)),\tag 2.11$$
where $\bold l'=\bold a_1^{-1}\bold la_1\in N^{-}$ and $P_1$ is $
P$ followed by the 
projection to $\frak u\ominus i\frak a_0$.  

Given $x\in \frak g$, if $x=x_{-}+x_{\frak h}+x_{+}$ is its triangular 
decomposition, then 
$$P(x)=-x_{+}^{*}+\frac 12(x_{\frak h}-x_{\frak h}^{*})+x_{+}.\tag 2.12$$
If $x\in \frak g_0$, then $x_{-}=x_{+}^{\sigma}$, and $x_{\frak h}
=x_{\frak t_0}\in \frak t_0$.  Because $\bold l'\bold a_1^{-1}$ 
maps $\frak n_{-}$ into itself, $(2.11)$ is given by 
$$P_1(Ad(\bold l'\bold a_1^{-1})(x))=-[(x_{+}^{\bold l'\bold a_1^{
-1}})_{+}]^{*}+(x_{\frak t_0}+(x_{+}^{\bold l'\bold a_1^{-1}})_{\frak t_
0})+(x_{+}^{\bold l'\bold a_1^{-1}})_{+}.\tag 2.13$$
Thus the determinant of $(2.10)$ is the same as the (real) 
determinant of the map $x_{+}\to (x_{+}^{\bold l'\bold a_1^{-1}})_{
+}$.  Because of the 
unipotence of $Ad(\bold l')$, this is equal to 
$$\prod_{\alpha >0}\vert \bold a_1^{-\alpha}\vert^2=\bold a_1^{-4
\delta}=a_{\phi}^{2\delta}.\tag 2.14$$
\qed
\enddemo

We will now show that the integral $(2.8)$ can be 
computed using a Duistermaat-Heckman exact stationary 
phase calculation.  The relevant Poisson structure was 
discovered by Evens and Lu in a very general setting 
([EL]).  We will introduce this structure directly, but to 
understand why it is natural the reader will need to 
consult the original paper.  

To do calculations we will use the isomorphism of 
vector bundles 
$$G_0\times_K\frak p\to T(G_0/K):[g_0,x]\to\frac d{dt}\vert_{t=0}
(g_0e^{tx}K),$$
and we will use the Killing form to identify $\frak p^{*}$ with $
\frak p$.  

Consider the $Ad(T_0)$ and $Ad(A_0)$-stable decomposition of $\frak g$ 
as a direct sum of subalgebras:  
$$\frak g=\frak g_0\oplus (\frak n^{-}\oplus i\frak h_0).\tag 2.15$$
Let $pr_{\frak g_0}$ denote the projection $\frak g\to \frak g_0$ along this 
decomposition.  Given $x\in \frak g$, with triangular 
decomposition $x=x_{-}+x_0+x_{+}$, 
$$pr_{\frak g_0}(x)=(x_{+}^{\sigma}+(x_0)_{\frak h_0}+x_{+}).\tag 2.16$$

The Evens-Lu Poisson bivector is given by 
$$\Pi ([g_0,x]\wedge [g_0,y])=\langle\Omega (g_0)(x),y\rangle ,\tag 2.17$$
where $\Omega (g_0):\frak p\to \frak p$ is given by 
$$\Omega (g_0)(x)=\{(pr_{\frak g_0}(ix^{g_0}))^{g_0^{-1}}\}_{\frak p}
.\tag 2.18$$
The operator $\Omega$ satisfies the equivariance condition 
$$\Omega (a_0g_0k)=Ad(k)^{-1}\Omega (g_0)Ad(k).\tag 2.19$$

To understand $\Omega$, it is useful to consider the augmented 
operator $\tilde{\Omega }:\frak g_0\to \frak g_0$ given by 
$$\tilde{\Omega }(g_0)(x_{\frak k}+x_{\frak p})=\{(pr_{\frak g_0}
((x_{\frak k}+ix_{\frak p})^{g_0}))^{g_0^{-1}}\}.\tag 2.20$$
Relative to the decomposition $\frak g_0=\frak k\oplus \frak p$, 
$$\tilde{\Omega }=\left(\matrix 1&*\\
0&\Omega\endmatrix \right).\tag 2.21$$
This augmented operator can be factored as the 
composition of four operators 
$$\frak g_0@>{I}>>\frak u@>{Ad(\bold u(g_0))}>>\frak u@>{T}>>\frak g_
0@>{Ad(\bold a_0^{-1}g_0)^{-1}}>>\frak g_0\tag 2.22$$
where the first operator is given by $I(x_{\frak k}+x_{\frak p})=
x_{\frak k}+ix_{\frak p}$, 
and $T(g_0)$ maps $x=-x^{*}_{+}+(x_{t_0}+x_{i\frak a_0})+x_{+}$ to 
$$T(g_0)(x)=pr_{\frak g_0}(x^{\bold l'\bold a_1(g_0)})=$$
$$[(x_{+}^{\bold l'\bold a_1})_{+}]^{\sigma}+(x_{t_0}+(x^{\bold l'
\bold a_1}_{+})_{\frak t_0+\frak a_0})+(x_{+}^{\bold l'\bold a_1}
)_{+},\tag 2.23$$
where $\bold l'=\bold a_0\bold l\bold a_0^{-1}$, and the last equality depends upon the 
fact that conjugation by $\bold l'\bold a_1(g_0)$ maps $\frak n^{
-}$ into itself, and 
that $\frak n^{-}$ terms disappear when we use $(2.16)$.  
 
\proclaim{ (2.24) Lemma} (a) $\Omega\in so(\frak p)$; the Schouten 
bracket $[\Pi ,\Pi ]$ vanishes, so that $(G_0/K,\Pi )$ is a Poisson 
manifold.
  
(b) $ker(\Omega (g_0))=\{[g_0,(\frak a_0^{\bold u(g_0)^{-1}})_{\frak p}
]\}$ 

(c) $Pfaffian(\Omega (g_0)\vert_{ker(\Omega )^{\perp}})=\bold a_1
(g_0)^{2\delta}$.  
\endproclaim

\demo{Proof of (2.24)} For (a) let $X=x^{g_0}$, $Y=y^{g_0}$, 
$x,y\in \frak p$.  Then 
$$\langle\Omega (g_0)x,y\rangle =\langle pr_{\frak g_0}(iX),Y\rangle$$
$$=\langle -iX_{+}^{\sigma}+iX_{+},Y_{+}^{\sigma}+Y_{\frak h_0}+Y_{
+}\rangle =2\langle iX_{+},Y_{+}^{\sigma}\rangle .\tag 2.25$$
This is clearly skew-symmetric in $X$ and $Y$, because $\sigma$ 
preserves the Killing form and it is complex antilinear.  
For the second part of $(a)$ we refer to [EL] (or see \S 3 
of [FO] for an exposition specific to this case).

For (b), note that $(2.23)$ implies the kernel of $T$ is $i\frak a_
0$.  
Thus $(2.22)$ implies 
$$ker(\tilde{\Omega }(g_0))=\{[g_0,x]:(x_{\frak k}+ix_{\frak p})\in 
i\frak a_0^{\bold u(g_0)^{-1}}\}\tag 2.26$$
This, together with $(2.21)$, implies $(b)$.  

For (c), note that in $(2.22)$ the first, second and fourth 
operators preserve volume determined by the Killing 
form.  The determinant of $T$ (relative to the Killing 
form volumes) is the same as the determinant of the 
operator on $\frak n^{+}$ which maps $x_{+}$ to $(x_{+}^{\bold l'
\bold a_1})_{+}$.  This 
determinant equals 
$$\prod_{\alpha >0}\bold a_1^{2\alpha}=\bold a_1^{4\delta}\tag 2.27$$
Thus the Pfaffian is $\bold a_1^{2\delta}$.  
\qed
\enddemo

By $(b)$ the tangent directions in $G_0/K$ determining the 
symplectic leaves are given by $[g_0,x]$ such that $x^{\bold u}\perp 
\frak a_0$.  
This is clearly $A_0$-invariant, because $\bold u(a_0g_0)=\bold u
(g_0)$.  
Thus the left action of $A_0$ permutes the symplectic 
leaves.  The symplectic form is given by the formula 
$$\omega ([g_0,x],[g_0,y])=\langle\Omega (g_0)\vert_{ker(\Omega )^{
\perp}})^{-1}(x),y\rangle .\tag 2.28$$
This form does not in general descend to a form on
the quotient $A_0\backslash G_0/K$.  However $(c)$ of the preceding 
Lemma does imply that the volume form descends.

\proclaim{ (2.29) Proposition}(a) The action of $T_0$ is 
Hamiltonian with momentum map 
$$\mu :G_0/K\to (\frak t_0)^{*}:g_0K\to\langle ilog(\bold a_1(g_0
)),\cdot\rangle ,$$
This momentum map is proper, and it is semibounded.  

(b) The symplectic measure is 
$$\omega^d/d!=\bold a_1(g_0K)^{-2\delta}dV_{A_0\backslash G_0/K}(
A_0g_0K)$$
(where the invariant measure is suitably normalized).  

\endproclaim

\demo{Proof of (2.29)}Part (a) is proven in [FO] (Lemma 
$3.3$,
which in turn refers to a result of Van Den Ban).  
Part (b) follows from (c) of 
$(2.24)$.  \qed
\enddemo

We can now apply the Duistermaat-Heckman exact 
stationary phase method, as generalized to 
noncompact manifolds 
in [PW]. For definiteness we will consider the 
symplectic leaf through the basepoint of $G_0/K$.

We must first find the fixed points of the $T_0$ 
action.  Suppose that $g_0K$ is fixed by $T_0$. If 
we choose $\lambda\in T_0$ 
which generates a dense subgroup of $T_0$, this is 
equivalent to 
$g_0^{-1}\lambda g_0\in K$. Since $T_0$ is maximal 
abelian in $K$, we can assume (by multiplying $g_0$ on the 
right by $k\in K$ if necessary) that $g_0^{-1}\lambda g_0\in T_0$. Since 
$N_{G_0}(T_0)=N_K(T_0)exp(\frak h_0)$, $g_0K=a_0K$ for some $a_0\in 
A_0$. 
Thus each symplectic leaf has exactly one $T_0$ fixed 
point. Since we are considering the leaf through the 
basepoint,
there is just one $T_0$ fixed point, the basepoint. 

If $X$ denotes the element of $\frak t_0$ corresponding to $\delta 
+\Lambda$ 
($\Lambda =-i\lambda$), then 
the Pfaffian of the infinitesimal action of $X$ at the 
basepoint equals 
$$Pf(ad(X)\vert_{\frak p})=\prod\langle\delta +\Lambda ,\alpha\rangle 
,\tag 2.30$$
where the product is over pairs of positive roots 
$(\alpha ,\Theta (\alpha ))$ which are not of compact type.  The 
Duistermaat-Heckman formula now implies $(2.3)$ in the 
case $\bold w=1$. This concludes 
the proof of $(2.2)$.  

\bigskip

We end this section with two brief remarks.
First, it is interesting to consider the integral 
$$\psi_{\Lambda}(g)=\int_{A_0\backslash G_0}\bold a_1(g_0g)^{-2\delta 
-2(\Lambda +\delta )}dg_0,\tag 2.31$$
for $g\in G$, which generalizes $(2.8)$.  
When this is well-defined, (1) this is a 
function of $g\in G_0\backslash G/U$, (2) this is a 
$G_0$-invariant eigenfunction for $G$-invariant 
differential operators on 
$G/U$; see Lemma 5.15 of ch2 section 5 of [H2] (one is 
usually interested in $U$-invariant eigenfunctions, i.e. 
spherical functions).  

To explicitly evaluate $(2.31)$, first note 
that $G_0exp(i\frak t_0)U=G$ (this is existence of polar 
decomposition for the nonRiemannian symmetric space
$G_0\backslash G$).  Thus we can suppose 
that $g=a\in exp(i\frak t_0)$.  In this case $(2.31)$ is an integral 
over $A_0\backslash G_0/C_K(a)$.  One can define a Poisson structure 
on $G_0/C_K(a)$, using the Evens-Lu method, as above 
(see \S 3 of [FO]). As in $(2.29)$, the momentum map can be 
identified with $log(\bold a_1(ga))$, and the symplectic volume of 
a symplectic leaf
can be identifed with the form
$\bold a_1(g_0a)^{-2\delta}dV_{A_0\backslash G_0/C_K(a)}$, via the projection 
to the double coset space. The fixed 
points for the 
$T_0$ action are of the form $a_0wC_K(a)$, where $a_0\in A_0$ and 
$w\in W(K,T_0)/W(C_K(a),T_0)$ (see Prop 4.3 of [FO]). Thus $(2.31
)$ 
equals
$$\sum_{\{w\}}\frac {\bold a_1(wa)^{-2(\delta +\Lambda )}}{\prod^
w\langle (\delta +\Lambda )^{w^{-1}},\alpha\rangle}$$
where given $w\in W(K,T_0)/W(C_K(a),T_0)$, the product is 
over (1) pairs of positive roots $(\alpha ,\Theta (\alpha ))$ which are not 
of compact type (relative to $\Theta$), and (2) positive compact 
type roots which vanish on $C_{\frak k}(a)$.

The second remark is that there is a kind of ``dual'' 
Poisson 
structure, on all of $U/K$, which can be used so that the 
sum $(2.4)$ has the structure of an exact stationary 
phase calculation. In the terminology of the paper [EL],
in this section we used the Lagrangian splitting 
$(2.15)$, to obtain a Poisson structure on $G_0/K$; 
the ``dual'' is the (Iwasawa) Lagrangian 
splitting $\frak g=\frak u\oplus (\frak h_{\Bbb R}+\frak n^{-})$, which induces a Poisson 
structure on $U/K$. This will hopefully be taken up 
elsewhere.

\bigskip

\centerline{Appendix. Special Features of the Group Case.}

\bigskip

In this appendix we will present a proof of $(0.1)$, using 
facts about Bott-Samelson resolutions of Schubert 
varieties.  One rationale for including this appendix is 
that many of the arguments are valid in the more 
general context of Kac-Moody Lie algebras and groups.  
Throughout this appendix, we will use the notation and 
basic results in [Kac].  

We start with the following data:  $A$ is an irreducible 
symmetrizable generalized Cartan matrix; $\frak g=\frak g(A)$ is the 
corresponding Kac-Moody Lie algebra, realized via its 
standard (Chevalley-Serre) presentation; $\frak g=\frak n^{-}\oplus 
\frak h\oplus \frak n^{+}$ is 
the triangular decomposition; $\frak b=\frak h\oplus \frak n^{+}$ the upper Borel 
subalgebra; $G=G(A)$ is the algebraic group associated to 
$A$ by Kac-Peterson; $H,N^{\pm}$ and $B$ are the subgroups of $G$ 
corresponding to $\frak h$, $\frak n^{\pm}$, and $\frak b$, respectively; $
K$ is the 
``unitary form" of $G$; $T=K\cap H$ the maximal torus; and 
$W=N_K(T)/T\simeq N_G(H)/H$ is the Weyl group.  

A basic fact is that $(G,B,N_G(H))$ with Weyl group $W$ is 
an abstract Tits system.  This yields a complete 
determination of all the (parabolic) subgroups between $B$ 
and $G$.  They are described as follows.  

Let $\Phi$ be a fixed subset of the simple roots.  The 
subgroup of $W$ generated by the simple reflections 
corresponding to roots in $\Phi$ will be denoted by $W(\Phi )$.  
The parabolic subgroup corresponding to $\Phi$, $P=P(\Phi )$, is 
given by $P=BW(\Phi )B$.  

Given $w\in N_K(T)$, we will denote its image in $W/W(\Phi )$ by 
$\bar {w}$.  

The basic structural features of $G/P$ which we will 
need are the Birkhoff and Bruhat decompositions 
$$G/P=\bigsqcup\Sigma_{\bar {w}},\qquad\Sigma_{\bar {w}}=N^{-}wP\tag A.1$$
$$G/P=\bigsqcup C_{\bar {w}},\qquad C_{\bar {w}}=BwP,\tag A.2$$
respectively, where the indexing set is $W/W(\Phi )$ in both 
cases.  The strata $\Sigma_{\bar {w}}$ are infinite dimensional if $
\frak g$ is 
infinite dimensional, while the cells $C_{\bar {w}}$ are always finite 
dimensional.  Our main interest lies in the Schubert 
variety $\bar {C}_{\bar {w}}$, the closure of the cell.  

Fix $\bar {w}\in W/W(\Phi )$.  We choose a representative $w\in N
(T)$ 
of minimal length $n$; for definiteness we will always 
take $w$ of the form 
$$w=r_n\cdots r_1\tag A.3$$
where $r_i=i_{\alpha_i}(\left(\matrix 0&1\\
-1&0\endmatrix \right))$, and $i_{\alpha_i}:SL_2\rightarrow G$ is the canonical 
homomorphism of $SL_2$ onto the root subgroup 
corresponding to the simple root $\alpha_i$.  

\proclaim{ (A.4) Proposition}For $w$ as in (A.3), the map 
$$r_nexp(\frak g_{-\alpha_n})\times ..\times r_1exp(\frak g_{-\alpha_
1})\to G/P:(p_j)\to p_n..p_1P$$
is a complex analytic isomorphism onto $C_{\bar {w}}$.  
\endproclaim

This result is essentially (5) of [Kac] together with 
Tits's theory.  We will sketch a proof for completeness.  

\demo{Proof of (A.4)} Let $\Delta^{+}$ denote the positive roots, 
$\Delta^{+}(\Phi )$ the positive roots which are combinations of 
elements from $\Phi$.  The ``Lie algebra of $P$'' is $\frak p=\Sigma 
\frak g_{-\beta}\oplus \frak b$ 
where the sum is over $\beta\in\Delta^{+}(\Phi )$; this is the Lie algebra 
of $P$ in the sense that it is the subalgebra generated by 
the root spaces $\frak g_{\gamma}$ for which $exp:\frak g_{\gamma}
\rightarrow G$ is defined and 
have image contained in $P$.  The subgroups $exp(\frak g_{\gamma}
)$ 
generate $P$.  We also let $\frak p^{-}$ denote the subalgebra 
opposite $\frak p$:  $\frak p=\sum \frak g_{-\gamma}$, where the sum is over 
$\gamma\in\Delta^{+}\setminus\Delta^{+}(\Phi )$.  The corresponding group will be denoted 
by $P^{-}$.  

The cell $C_w$ is the image of the map $N^{+}\rightarrow G/P:u\rightarrow 
uwP$.  
The stability subgroup at $wP$ is $N^{+}\cap wPw^{-1}$.  

At the Lie algebra level we have the splitting 
$$\frak n^{+}=\frak n^{+}\cap Ad(w)(\frak p)\oplus \frak n^{+}\cap 
Ad(w)(\frak p^{-}).\tag A.5$$
The second summand equals 
$$\frak n_w^{+}=\oplus\,\frak g_{\beta}\tag A.6$$
where the sum is over roots $\beta >0$ with 
$w^{-1}\beta\in -(\Delta^{+}\setminus\Delta^{+}(\Phi ))$.  These roots $
\beta$ are necessarily real, 
so that $exp:\frak n_w^{+}\rightarrow N_w^{+}\subseteq N^{+}$ is well-defined.  

For $q\in \Bbb Z^{+}$ let $N_q^{+}$ denote the subgroup corresponding to 
$\frak n_q^{+}=span\{\frak g_{\beta}:height(\beta )\geq q\}$.  Then $
N^{+}/N_q^{+}$ is a finite 
dimensional nilpotent Lie group, and it is also simply 
connected.  By taking $q$ sufficiently large and considering 
the splitting (A.5) modulo $\frak n_q^{+}$, we conclude by finite 
dimensional considerations that each element in $N^{+}$ has a 
unique factorization $n=n_1n_2$, where $n_1\in N_w^{+}$ and 
$n_2\in N^{+}\cap wPw^{-1}$:  
$$N^{+}\simeq N_w^{+}\times (n^{+}\cap wPw^{-1}).\tag A.7$$

The important point here is that modulo $N_q^{+}$ we can 
control $N^{+}\cap wPw^{-1}$ by the exponential map.  

We now recall the following standard 

\proclaim{ (A.8) Lemma} In terms of the minimal 
factorization $w=r_n\cdots r_1$, the roots $\beta >0$ with $w^{-1}
\beta <0$ 
are given by 
$$\beta_j=r_n\cdots r_{j+1}(\alpha_j)=r_n\cdots r_j(-\alpha_j),\qquad 
1\leq j\leq n.$$
\endproclaim

Because $w$ is a representative of $\bar {w}\in W/W(\Phi )$ of minimal 
length, all of these $\beta_j$ satisfy $w^{-1}\beta_j\in -(\Delta^{
+}\setminus\Delta^{+}(\Phi ))$.  
Otherwise, if say $w^{-1}\beta_j\in -\Delta^{+}(\Phi )$, then 
$$w^{-1}r_{\beta_j}w=r_1\cdots r_{j-1}r_jr_{j-1}\cdots r_1\in N(T
)\cap P\tag A.9$$
and $w'=w(w^{-1}r_{\beta_j}w)=r_n\cdots\hat {r}_j\cdots r_1$ would be a 
representative of $\bar {w}$ of length $<n$ (here we have used 
the fact that $W(\Phi )=N(T)\cap P/T$, which follows from the 
Bruhat decomposition).  For future reference we note 
this proves that 
$$N_w^{+}=N^{+}\cap (N^{-})^w=N^{+}\cap (P^{-})^w\tag A.10$$
and (2.4) shows that 
$$N_w^{+}\times w\cong C_{\bar {w}}.\tag A.11$$
Now for any $1\le p\le q\le n$, $\bigoplus_{p\le j\le p}g_{\beta_
j}$ is a subalgebra of 
$n_w^{+}$.  Thus by (2.7) 
$$exp(\frak g_{\beta_n})\times\cdots\times exp(\frak g_{\beta_1})
\times w\cong C_{\bar {w}}.\tag A.12$$
This yields (A.4) when we write 
$$exp(\frak g_{\beta_j})=r_n\cdots r_jexp(\frak g_{\alpha_j})r_j\cdots 
r_n.\tag A.13$$
\qed
\enddemo

We now note several important corollaries of (A.4).  

For each $i$, let $P_i$ denote the parabolic subgroup 
$i_{\alpha_i}(SL_2)B$.  Let 
$$\Gamma_w=P_n\times_B\cdots\times_BP_1/B\tag A.14$$
where 
$$P_n\times\cdots\times P_1\times B\times\cdots\times B\rightarrow 
P_n\times\cdots\times P_1\tag A.15$$
is given by 
$$(p_j)\times (b_j)\rightarrow (p_nb_n,b_n^{-1}p_{n-1}b_{n-1},\cdots 
,b_2^{-1}p_1b_1).\tag A.16$$
We have written ``$\Gamma_w$" instead of ``$\Gamma_{\bar {w}}$" to indicate that 
this compact complex manifold depends upon the 
factorization.  

\proclaim{ (A.17) Corollary}The map 
$$\Gamma_w\rightarrow\bar {C}_{\bar {w}}:(p_j)\rightarrow p_n\cdots 
p_1P$$
is a desingularization of $\bar {C}_{\bar {w}}$.  
\endproclaim

Let 
$$SL_2'=\{g=\left(\matrix a&b\\
c&d\endmatrix \right)\in SL(2,\Bbb C):a\ne 0\}.\tag A.17$$

\proclaim{ (A.18) Corollary} Let $\phi$ denote the surjective 
map 
$$SL_2\times\cdots\times SL_2\rightarrow\bar {C}_{\bar {w}}:(g_j)
\rightarrow r_ni_{\alpha_n}(g_n)\cdots r_1i_{\alpha_1}(g_1)P.$$
The inverse image of $C_{\bar {w}}$ under $\phi$ is $SL_2'\times\cdots
\times SL_2'$.  
\endproclaim

\demo{Proof of (A.18)} Let $\sigma =r_{n-1}\cdots r_1$.  It suffices to 
show that for the natural actions 
$$r_ni_{\alpha_n}(SL_2')\times C_{\bar{\sigma}}\rightarrow C_{\bar {
w}},\
\tag A.19$$
$$r_ni_{\alpha_n}(SL_2\setminus SL_2')\times C_{\bar{\sigma}}\rightarrow
\bar {C}_{\bar {w}},\tag A.20$$
and 
$$r_ni_{\alpha_n}(SL_2)\times (\bar {C}_{\bar{\sigma}}\setminus C_{
\bar{\sigma}})\rightarrow\bar {C}_{\bar {w}}\setminus C_{\bar {w}}
.\tag A.21$$
The first line, $(A.19)$, follows from (A.4) since 
$i_{\alpha_n}(SL_2')\subseteq exp(-\frak g_{-\alpha_n})B$ and $B\times 
C_{\bar{\sigma}}\subseteq C_{\bar{\sigma}}$.  The second line 
follows from 
$$r_ni_{\alpha_n}\left(\matrix 0&b\\
c&d\endmatrix \right)\cdot C_{\bar{\sigma}}=i_{\alpha_n}\left(\matrix 
c&b\\
0&d\endmatrix \right)\cdot C_{\bar{\sigma}}\subseteq C_{\bar{\sigma}}
.\tag A.22$$
For the third line it's clear that the image of the left 
hand side is a union of cells, since we can replace 
$r_ni_{\alpha_n}(SL_2)$ by $P_n$.  This image is at most $n-1$ 
dimensional.  Therefore it must have null intersection 
with $C_{\bar {w}}$.  \qed
\enddemo

Fix an integral functional $\lambda\in \frak h^{*}$ which is antidominant.  
Denote the (algebraic) lowest weight module 
corresponding to $\lambda$ by $L(\lambda )$, and a lowest weight vector 
by $\sigma_{\lambda}$.  Let $\Phi$ denote the simple roots $\alpha$ for which 
$\lambda (h_{\alpha})=0$, where $h_{\alpha}$ is the coroot, $P=P(
\Phi )$ the 
corresponding parabolic subgroup.  The Borel-Weil 
theorem in this context realizes $L(\lambda )$ as the space of 
strongly regular functions on $G$ satisfying 
$$f(gp)=f(g)\lambda (p)^{-1}\tag A.23$$
for all $g\in G$ and $p\in P$, where we have implicitly 
identified $\lambda$ with the character of $P$ given by 
$$\lambda (u_1w\exp(x)u_2)=\exp\lambda (x)\tag A.24$$
for $x\in \frak h,\;u_1,u_2\in N^{+},\;w\in W(\Phi )$.  Thus we can view $
L(\lambda )$ 
as a space of sections of the line bundle 
$$\Cal L_{\lambda}=G\times_{\lambda}\Bbb C\rightarrow G/P.\tag A.25$$
If $\frak g$ is of finite type, then $L(\lambda )=H^0(\Cal L_{\lambda}
)$; if $\frak g$ is affine 
(and untwisted), then $L(\lambda )$ consists of the holomorphic 
sections of finite energy, as in [PS].  

Normalize $\sigma_{\lambda}$ by $\sigma_{\lambda}(1)=1$.  

\proclaim{ (A.26) Proposition} Let $\bar {w}\in W/W(\Phi )$, and let 
$w=r_n\cdots r_1$ be a representative of minimal length $n$.  The 
positive roots mapped to negative roots by $w$ are given 
by 
$$\tau_j=r_1\cdots r_{j-1}(\alpha_j),\qquad 1\leq j\leq n;$$
let $\lambda_j=-\lambda (h_{\tau_j})$, where $h_{\tau}$ is the coroot corresponding to 
$\tau$.  Then 
$$\sigma_{\lambda}^w(r_ni_{\alpha_n}(g_n)\cdots r_1i_{\alpha_1}(g_
1))=\prod_1^na_j^{\lambda_j}$$
where $g=\left(\matrix a&b\\
c&d\endmatrix \right)\in SL_2$.  
\endproclaim

\demo{Proof of (A.26)} The claim about the $\tau_j$ is easily 
derived from (2.5).  None of these roots lie in $\Delta^{+}(\Phi 
)$, by 
the same argument as follows (2.5).  Thus each $\lambda_j>0$.  
It follows that $\Pi a_j^{\lambda_j}$ is nonzero precisely on the set 
$SL_2'\times\cdots\times SL_2'$.  

Now $\sigma_{\lambda}^w$, viewed as a section of $\Cal L_{\lambda}
\rightarrow G/P$, is nonzero 
precisely on the $w$-translate of the largest strata, 
$$w\Sigma_0=wP^{-}P=(P^{-})^wwP.\tag A.27$$
We claim the intersection of this with $\bar {C}_{\bar {w}}$ is $
C_{\bar {w}}$.  In 
one direction 
$$C_{\bar {w}}=\left(N^{+}\cap (P^{-})^w\right)wP\subseteq (P^{-}
)^wwP\tag A.28$$
by (2.6).  On the other hand $(N^{+}\cap (P^{-})^w)$ is a closed 
finite dimensional subgroup of $(P^{-})^w$.  Since $(P^{-})^w$ is 
topologically equivalent to $w\Sigma_0$, the limit points of $C_{
\bar {w}}$ 
must be in the complement of $w\Sigma_0$.  This establishes the 
other direction.  

It now follows from (A.4) that $\sigma_{\lambda}^w$ is also nonzero 
precisely on $SL_2'\times\cdots\times SL_2'$, viewed as a function of 
$(g_n,\cdots ,g_1)$.  

We now calculate that 
$$\sigma_{\lambda}^w(r_ni_{\alpha_n}(g_n)\cdots r_1i_{\alpha_1}(g_
1))=\sigma_{\lambda}(w^{-1}r_ni_{\alpha_n}(g_n)\cdots r_1i_{\alpha_
1}(g_1))\ $$
$$=\sigma_{\lambda}\left(\omega_{n-1}i_{\alpha_n}(g_n)\omega_{n-1}^{
-1}\omega_{n-2}i_{\alpha_{n-1}}(g_{n-1})\omega_{n-2}^{-1}\cdots\omega_
0i_{\alpha_1}(g_1)\omega_0^{-1}\right)\ $$
$$=\sigma_{\lambda}\left(i_{\tau_n}(g_n)i_{\tau_{n-1}}(g_{n-1})\cdots 
i_{\tau_1}(g_1)\right)\tag A.29$$
where we have set $\omega_i=r_1..r_i,\;0\leq i<n$, and we have 
used 
$$\omega_{i-1}(\alpha_i)=r_1\cdots r_{i-1}(\alpha_i)>0\tag A.30$$
to conclude that $\omega_{i-1}i_{\alpha_i}(g)\omega_{i-1}^{-1}=i_{
\tau_i}(g)$.  

The map 
$$SL_2\times\cdots\times SL_2\rightarrow w^{-1}\bar {C}_{\bar {w}}
:(g_j)\rightarrow i_{\tau_n}(g_n)\cdots i_{\tau_1}(g_1)P\tag A.31$$
is surjective and the inverse image of $\Sigma_0\cap w^{-1}\bar {
C}_{\bar {w}}$ is 
precisely $SL_2'\times\cdots\times SL_2'$.  

For $g=\left(\matrix a&b\\
c&d\endmatrix \right)\in SL_2'$, write $g=LDU$, where 
$$L=\left(\matrix 1&0\\
ca^{-1}&1\endmatrix \right),\quad D=\left(\matrix a&0\\
0&a^{-1}\endmatrix \right),\quad U=\left(\matrix 1&a^{-1}b\\
0&1\endmatrix \right).\tag A.32$$
Then for $(g_j)\in SL_2'\times\cdots\times SL_2'$, (3.2) equals 
$$\sigma_{\lambda}(i_{r_n}(L_nD_nU_n)\cdots i_{\tau_1}(L_1D_1U_1)
)\ $$
$$=\sigma_{\lambda}(i_{\tau_n}(L_nU_n')i_{\tau_{n-1}}(L_{n-1}'U_{
n-1}')\cdots i_{\tau_1}(L_1'U_1')i_{\tau_n}(D_n)\cdots i_{\tau_1}
(D_1))\
$$
$$=\sigma_{\lambda}(i_{r_n}(L_nU_n')\cdots i_{\tau_1}(L_1'U_1'))\Pi 
a_j^{\lambda_j}\tag A.33$$
where each $L_j'$ ($U_j'$) has the same form as $L_j$ ($U_j$, 
respectively).  This follows from the fact that $H$ 
normalizes each $exp(\frak g_{\pm r})$.  

Now each $L_j'U_j'\in SL_2'$, so that $i_{\tau_n}(L_nU_n')\cdots 
i_{\tau_1}(L_1'U_1')$ is in 
$\Sigma_0$.  We now conclude that 
$$\sigma_{\lambda}(i_{\tau_n}(L_nU_n')\cdots i_{\tau_1}(L_1'U_1')
)=1,\tag A.34$$
by the fundamental theorem of algebra, since this is 
polynomial and never vanishes.  \qed
\enddemo

\proclaim{ (A.35) Proposition}Suppose that $\frak g$ is finite 
dimensional.  Given $g\in K$ such that $gT\in\Sigma_1$, we can write 
$g$ uniquely as $g=lmau$, where $l\in N^{-}$, $m\in T$, $a\in exp
(\frak h_{\Bbb R})$, 
and $u\in N^{+}$.  Then 
$$\int_Ka(g)^{-i\lambda}=\prod_{\alpha >0}\frac {\langle 2\delta 
,\alpha\rangle}{\langle 2\delta -i\lambda ,\alpha\rangle},$$
where the integral is with respect to the normalized 
Haar measure of $K$, and $2\delta$ denotes the sum of the 
positive complex roots.  
\endproclaim

\demo{Proof of $(A.35)$} Let $\{\Lambda_j\}$ denote the set of basic 
dominant integral functionals.  

We apply $(A.26)$ to $w=w_0$.  We write $w_0=r_n..r_1$ as in 
$(A.3)$.  Then 
$$gT=i_{\tau_n}(g_n)i_{\tau_{n-1}}(g_{n-1})\cdots i_{\tau_1}(g_1)
T\tag A.36$$
$$a(g)=\prod_{j=1}^l\vert\sigma_{\Lambda_j}(g)\vert^{h_j}=\prod_{
j=1}^l(\prod_{k=1}^n\vert a_k\vert^{\Lambda_j(h_{\tau_k})})^{h_j}
=\prod_{k=1}^n\vert a_k\vert^{h_{\tau_k}},\tag A.37$$
since the $\Lambda_j$ are dual to the $h_j$.  Therefore 
$$a(g)^{-i\lambda}=\prod_{k=1}^n\vert a_k\vert^{-i\lambda (h_{\tau_
k})}.\tag A.38$$

Also, in terms of the coordinates $a_k$, the invariant 
measure is given by 
$$a(g)^{2\delta}\prod_k\vert a_k\vert^{-2}dm(a_k),\tag A.39$$
up to a normalization factor.  

The roots $\tau_k$ range over all the positive complex roots.  
Thus by $(A.26)$, 
$$\int_Ka(g)^{-i\lambda}=\int_{\Sigma_1}a(gT)^{-i\lambda}$$
$$\Cal Z^{-1}\prod_{\alpha >0}\int_{SU(2)}\vert a\vert^{(2\delta 
-i\lambda )(h_{\alpha})}\vert a\vert^{-2}=\Cal Z^{-1}\prod_{\alpha 
>0}\int_0^1r^{(2\delta -i\lambda )(h_{\alpha})-1}dr$$
$$=\Cal Z^{-1}\prod_{\alpha >0}\frac 1{(2\delta -i\lambda )(h_{\alpha}
)}=\prod_{\alpha >0}\frac {\langle 2\delta ,\alpha\rangle}{\langle 
2\delta -i\lambda ,\alpha\rangle}.$$
\qed
\enddemo

This proof was given a Poisson-theoretic interpretation 
in [Lu].  

\bigskip

\centerline{References}

\bigskip

\flushpar[EL] S Evens and J-H Lu, On the variety of 
Lagrangian subalgebras, Ann Scient Ec Norm Sup 4 (34) 
(2001) 631-668.  

\flushpar[FO] P Foth and M Otto, A symplectic 
realization of Van Den Ban's convexity theorem, 
math.SG/0505063.  

\flushpar[H1] S Helgason, Differential Geometry, Lie 
Groups, and Symmetric Spaces, Academic Press (1978) 

\flushpar[H2] --------, Groups and Geometric Analysis, 
Academic Press (1984).  

\flushpar[Kac] V Kac, Constructing groups from infinite 
dimensional Lie algebras, Infinite Dimensional Groups 
with Applications, edited by V Kac, MSRI publication, 
Springer-Verlag (1985) 167-216.  

\flushpar[Kn] A Knapp, Lie groups beyond an introduction, 
2nd edition, Birkhauser (2002).  

\flushpar[Lu] J Lu, Coordinates on Schubert cells, 
Kostant's harmonic forms, and the Bruhat-Poisson 
structure on $G/B$, Transf.  Groups 4, No.  4 (1999) 
355-374.  

\flushpar[Pi1] D Pickrell, Invariant measures for unitary 
forms of Kac-Moody Lie groups, Memoirs of the AMS, 
Vol 146, No 693 (2000).  

\flushpar[Pi2] ---------, An invariant measure for the 
loop space of a simply connected compact symmetric 
space, submitted to J Funct Anal.  

\flushpar[PW] E Prato and S Wu, Duistermaat-Heckman 
measures in a noncompact setting, Compositio Math.
94 no. 2 (1994) 113-128.  

\flushpar[Wolf] J Wolf, The action of a real semisimple 
group on a complex flag manifold.  I:  Orbit structure 
and holomorphic arc components, Bull A.M.S.  75 (1969) 
1121-1237.  

\flushpar[WZ] J Wolf and R Zierau, Cayley transforms 
and orbit structure in complex flag manifolds, Transf.  
Groups 2 (1997) no.  4, 391-405.  

\end

A simple but important fact, which we will use 
repeatedly, is that for $\bold w\in N_U(T)$, the equality $\bold w^{
*\Theta}=\bold w$ 
implies that $Ad(\bold w)^{-1}\Theta$ ($=\Theta Ad(\bold w)$) and $
Ad(\bold w)\Theta$ induce 
involutions on $\frak g$, $\frak u$, and $\frak t$.  For example there will be a 
splitting 
$$\frak t=\{Ad(\bold w)\Theta\vert_{\frak t}=1\}\oplus \{Ad(\bold w
)\Theta\vert_{\frak t}=-1\}.\tag 1.10$$
This will induce a splitting of $T$, modulo $T^{(2)}$, elements 
of order two.  

$================================$

Suppose that $\bold w_1\in N_U(T_0)$.  We must show that 
$\bold w=\bold w_1\bold w_1^{-\Theta}\in T_0^{(2)}$.  Suppose that $
t\in T_0$.  Since $t^{\Theta}=t$, 
$$\bold wt\bold w^{-1}=\bold w_1(\bold w_1^{-1}t\bold w_1)^{\Theta}
\bold w_1^{-1}=t$$
because $\bold w_1^{-1}t\bold w_1\in T_0$ also.  This should force $
\bold w\in T_0$.  
Since $\bold w^{\Theta}=\bold w$ and $\bold w^{-\Theta}=\bold w$, this implies $
\bold w$ has order 
two.  

Suppose $T_0\cdot uK=uK$.  Then $u^{-1}T_0u\subset K$.  Choose a 
regular element in $T_0$.  We can then find $k\in K$ such that 
$(uk)^{-1}$ will conjugate this regular element into $T_0$.  This 
implies that $uK=\bold w_1K$ for some $\bold w_1\in N_U(T_0)$.  

$=========================$ 

In $\S 3$ a number of classical cases are considered.  It 
turns out that in Cayley coordinates the integrands are 
computable and rational.  This is of interest in 
connection with the large rank limit, which corresponds 
to understanding the analogue of $(0.4)$ for an infinite 
rank unitary type symmetric space.  This will hopefully 
be worked out in Derek Habermas's dissertation.  

$==========================================$ 

\centerline{\S 3. Classical Examples and }
\centerline{the Mysterious Role of Cayley Coordinates}

\bigskip

Consider the case of the $2$-sphere.  The Iwasawa 
factorization of $g_0\in SU(1,1)$ is given by 
$$g_0=\left(\matrix a&b\\
\bar {b}&\bar {a}\endmatrix \right)=\left(\matrix 1&0\\
l&1\endmatrix \right)\left(\matrix \alpha&0\\
0&1/\alpha\endmatrix \right)\left(\matrix a/\alpha&b/\alpha\\
-\bar {b}/\alpha&\bar {a}/\alpha\endmatrix \right)$$
where $\alpha =(\vert a\vert^2+\vert b\vert^2)^{1/2}$, $l=\frac {
\bar {b}}a(\frac {\alpha^2+1}{\alpha^2})$.  Thus 
$$\Phi (\left(\matrix a&b\\
\bar {b}&\bar {a}\endmatrix \right))=\phi (\left(\matrix a/\alpha&
b/\alpha\\
-\bar {b}/\alpha&\bar {a}/\alpha\endmatrix \right))=\frac 1{\alpha^
2}\left(\matrix 1&2ab\\
-2\bar {a}\bar {b}&1\endmatrix \right)$$
$$=\left(\matrix 1&0\\
\frac {-2\bar {a}\bar {b}}{\alpha^2}&1\endmatrix \right)\left(\matrix 
\alpha^{-2}&0\\
0&\alpha^2\endmatrix \right)\left(\matrix 1&\frac {2ab}{\alpha^2}\\
0&1\endmatrix \right).$$
In terms of $z=-\bar {b}/a$, which satisfies $\vert z\vert <1$, we have 
$$\alpha^{-2}=\frac {\vert a\vert^2-\vert b\vert^2}{\vert a\vert^
2+\vert b\vert^2}=\frac {1-\vert z\vert^2}{1+\vert z\vert^2}$$
and $(*)$ equals 
$$=\left(\matrix 1&0\\
\frac {2z}{1+\vert z\vert^2}&1\endmatrix \right)\left(\matrix \frac {
1-\vert z\vert^2}{1+\vert z\vert^2}&0\\
0&\frac {1+\vert z\vert^2}{1-\vert z\vert^2}\endmatrix \right)\left
(\matrix 1&\frac {-2\bar {z}}{1+\vert z\vert^2}\\
0&1\endmatrix \right).$$
Note that 
$$a_{\phi}=(\frac {1-\vert z\vert^2}{1+\vert z\vert^2})^{h_1},\quad 
v=\frac {1-\vert z\vert^2}{1+\vert z\vert^2}<1.$$

We write $\lambda =\lambda\alpha_1$.  Note $\delta =\delta'=\frac 
12\alpha_1$.  Then 
$$\int_{\Sigma_1^{\phi (SU(2)/U(1))}}a_{\phi}(g)^{-i\lambda}=\int_{
SU(1,1)/U(1)}\bold a^{2i\lambda}\Phi^{*}(\delta_{U/K})$$
$$=\int_{\vert z\vert <1}(\frac {1-\vert z\vert^2}{1+\vert z\vert^
2})^{-2i\lambda}(1+\vert z\vert^2)^{-2}dm(z)$$
$$=\int_{v=0}^1v^{-2i\lambda}dv=\frac 1{1-2i\lambda}v^{1-2i\lambda}
\vert_{v=0}^1=\frac {\langle\delta ,\alpha_1\rangle}{\langle\delta 
-i\lambda\alpha_1,\alpha_1\rangle}.$$

This calculation can be refined in the following way 
(although the details become more involved, and use very 
special features of this low rank case).  

We want to compare the integrals 
$$\int_{G_0}\bold a(g_0g)^{-2\delta -2(\delta +\Lambda )}dg_0\qquad 
vs\qquad\int_U\bold a(ug)^{-2\delta -2(\delta +\Lambda )}du$$
The first integral is a function on $G_0\backslash G/U$ or $A/W(K
)$, 
whereas the second integral is a function on $U\backslash G/U$ or 
$A/W$.  

Consider the case of $SU(1,1)$.  Let 
$$a_1=\left(\matrix A&0\\
0&A^{-1}\endmatrix \right)$$
The calculation above can be generalized to 
$$\int_{SU(1,1)}\bold a(g_0a_1)^{-2\delta -2(\delta +\Lambda )}=\frac {
\langle\delta ,\alpha_1\rangle}{\langle\delta +\Lambda ,\alpha_1\rangle}\frac {
a_1^{-2\Lambda +\alpha_1}}{a_1^{\alpha_1}+a_1^{-\alpha_1}}$$
whereas 
$$\int_{SU(2)}\bold a(ua_1)^{-2\delta -2(\delta +\Lambda )}du=\frac {
\langle\delta ,\alpha_1\rangle}{\langle\delta +\Lambda ,\alpha_1\rangle}\frac {
a_1^{-2\Lambda -\alpha_1}-a_1^{-2\Lambda +\alpha_1}}{a_1^{\alpha_
1}-a_1^{-\alpha_1}}$$
In the first case we get the Weyl dimension by simply 
evaluating at $a_1=1$, whereas in the second case it arises 
in the asymptotics at infinity.  

\bigskip

We now consider Cayley coordinates.   A remarkable fact 
is that there are quite elegant formulas for the diagonal 
in Cayley coordinates.  It is worth remarking that the 
very formulation of an invariant measure, and a diagonal 
distribution conjecture, in infinite rank, depends upon 
the existence of Cayley coordinates (see [Pi]).  So this is 
a very special choice of coordinates.  

We begin with the case of $X=Gr(n,\Bbb C^{n+m})$.  In this case 
$\frak g=sl(n+m,\Bbb C)$, $\Theta$ is conjugation by $diag(1_n,-1_
m)$, and we 
choose the usual triangular decomposition, with $\frak h$ 
consisting of diagonal matrices.  

We parameterize $\phi (U/K)$ using Cayley coordinates (this 
is equivalent to graph coordinates for the Grassmannian, 
which has a natural generalization to any Hermitian 
symmetric space) 
$$i\frak p\to\phi (U/K):X=\left(\matrix 0&-Z^{*}\\
Z&0\endmatrix \right)\to g=(1-X)(1+X)^{-1}$$
$$=\left(\matrix \frac {1-Z^{*}Z}{1+Z^{*}Z}&2Z^{*}(1+ZZ^{*})^{-1}\\
-2Z(1+Z^{*}Z)^{-1}&\frac {1-ZZ^{*}}{1+ZZ^{*}}\endmatrix \right)$$
$$=\left(\matrix 1&0\\
-2Z(1-Z^{*}Z)^{-1}&1\endmatrix \right)\left(\matrix \frac {1-Z^{*}
Z}{1+Z^{*}Z}&0\\
0&\frac {1+ZZ^{*}}{1-ZZ^{*}}\endmatrix \right)\left(\matrix 1&2(1
-Z^{*}Z)^{-1}Z^{*}\\
0&1\endmatrix \right)$$
We must compute the integral 
$$\frac 1{\Cal Z}\int a(\left(\matrix \frac {1-Z^{*}Z}{1+Z^{*}Z}&
0\\
0&\frac {1+ZZ^{*}}{1-ZZ^{*}}\endmatrix \right))^{-i\lambda}det(1+
Z^{*}Z)^{-n-m}dm(Z)$$

Suppose that $n=1$.  In this case $Z$ is a column vector 
with entries $z_i$.  A Maple calculation shows that 
$$a_{\phi}=diag(\frac {1-\sum_1^m}{1+\sum_1^m},\frac {1+\sum_1^1-
\sum_2^m}{1-\sum_1^m},\frac {1+\sum_1^2-\sum_3^m}{1+\sum_1^1-\sum_
2^m},..,\frac {1+\sum_1^m}{1+\sum_1^{m-1}-\sum_m^m})=$$
$$=(1-\sum_1^m)^{h_1}(1+\sum_1^1-\sum_2^m)^{h_2}..(1+\sum_1^{m-1}
-\sum_m^m)^{h_m}(1+\sum_1^m)^{-(h_1+..+h_m)}$$
where $\sum_j^k=\sum_j^k\vert z_i\vert^2$.  From this we see that the 
intersection with $\Sigma^G_1$ has $m+1$ connected components, 
corresponding to the following elements of $T^{(2)}$:  the 
identity and the $m$ elements 
$$\bold w=exp(i\pi h_{\alpha_1})exp(i\pi h_{\alpha_2})..exp(i\pi 
h_{\alpha_j})=diag(-1,1,..,1,-1,1,..,1)$$
where the second $-1$ occurs in the $j+1$ slot, $1\le j\le m$.  
Note that we can also write 
$$a_{\phi}=(\frac {1-\sum_1^m}{1+\sum_1^m})^{h_1}(\frac {1+\sum_1^
1-\sum_2^m}{1+\sum_1^m})^{h_2}..(\frac {1+\sum_1^{m-1}-\sum_m^m}{
1+\sum_1^m})^{h_m}=\prod v_j^{h_j}$$
 
In the case $\bold w=1$, $0<v_1\le v_2\le ..\le v_m<1$, which is a 
convex set.  Moreover we have 
$$dv_1\wedge ..\wedge dv_m=c\frac 1{(1+\sum_1^mu_j)^{1+m}}du_1\wedge 
..\wedge du_m$$
In the case $\bold w=1$, the integral equals (up to a 
normalization constant) 
$$\int_{0<v_1\le ..\le v_m<1}\prod_1^mv_j^{-i\lambda (h_j)}dv_1\wedge 
..\wedge dv_m.$$

We consider the component corresponding to $\bold w=1$.  Let 
$u_j=\vert z_j\vert^2$.  Modulo a normalization factor, the integral 
equals 
$$\int (\frac {1-\sum_1^m}{1+\sum_1^m})^{-i\lambda_1}(\frac {1+\sum_
1^1-\sum_2^m}{1-\sum_1^m})^{-i\lambda_2}..(\frac {1+\sum_1^m}{1+\sum_
1^{m-1}-\sum_m^m})^{-i\lambda_{m+1}}(1+\sum_1^m)^{-1-m}du_m..du_1$$
where the domain of integration is the standard 
$m$-simplex, $0\le u_i$, $\sum u_i\le 1$ (and the summations are now 
over the $u_j$).  This equals 
$$\int (1-\sum_1^m)^{-i\alpha_1}(1+\sum_1^1-\sum_2^m)^{-i\alpha_2}
..(1+\sum_1^{m-1}-\sum_m^m)^{-i\alpha_m}(1+\sum_1^m)^{-1-m+i\sum\alpha_
j}$$
where we have written $\alpha_i=\lambda_i-\lambda_{i+1}$.  We claim that 
this integral equals 
$$=\frac 1{(1-i\alpha_1)(2-i(\alpha_1+\alpha_2)..(m-i(\alpha_1+..
+\alpha_m))}.$$

To inductively evaluate this integral, let $s=u_{m-1}+u_m$ 
and $t=u_{m-1}-u_m$.  In terms of these variables, the new 
constraints are $\sum_1^{m-2}+s\le 1$ and $-s\le t\le s$.  When we 
substitute these variables into $()$, $t$ appears only in the 
term with exponent $-i\alpha_m$.  We then first integrate with 
respect to $t$.  Then $()$ equals 
$$\int_{\{u_i:i\le m-2\},s}(1-\sum_1^{m-2}-s)^{-i\alpha_1}..(1+\sum_
1^{m-2}-s)^{-i\alpha_{m-1}}$$
$$\{\frac 1{1-i\alpha_m}(1+\sum_1^{m-2}+t)^{1-i\alpha_m}\vert_{t=
-s}^{t=s}\}(1+\sum_1^{m-2}+s)^{-1-m+i\sum\alpha_j}$$
with constraints $0\le u_i,s$ and $\sum u_i+s\le 1$.  This equals 
$(1-i\alpha_m)$ times the difference of two integrals which we 
can evaluate using the induction hypothesis.  We obtain 
$$\frac 1{1-i\alpha_m}\{\frac 1{(1-i\alpha_1)(2-i(\alpha_1+\alpha_
2))..(m-1-i(\alpha_1+..+\alpha_{m-1}))}$$
$$-\frac 1{(1-i\alpha_1)..(m-2-i(\alpha_1+..+\alpha_{m-2})(m-1-(-
1+i(\alpha_1+..+\alpha_m))}\}.$$
After obtaining a common denominator, we obtain the 
claim.  

We can also calculate that 

$$\prod_1^m\frac {\langle\delta ,\lambda_1-\lambda_{1+j}\rangle}{
\langle\delta -i\lambda ,\lambda_1-\lambda_{1+j}\rangle}=\frac {m
!}{(1-i(\lambda_1-\lambda_2))..(m-i(\lambda_1-\lambda_{m+1}))}$$
$$=\frac {m!}{(1-i\alpha_1)(2-i(\alpha_1+\alpha_2))..(m-i(\alpha_
1+..+\alpha_m))}$$
When we take account of the normalization, we obtain a 
proof of the conjecture in the case of $\Bbb C\Bbb P^m$.  

In the case of $Gr(2,\Bbb C^4)$, we have 
$$a_{\phi}=(1-u_{31}+u_{32}-u_{41}+u_{42}-\vert detZ\vert^2)^{h_1}
(1-u_{31}-u_{32}-u_{41}-u_{42}+\vert detZ\vert^2)^{h_2}$$
$$*(1+u_{31}+u_{32}-u_{41}-u_{42}-\vert detZ\vert^2)^{h_3}det(1+Z^{
*}Z)^{-h_1-h_2-h_3}$$
The integration in this case is not so easy.  

We now consider the case of $S^{2n}$.  We realize 
$\frak g=\frak o(2n+1,\Bbb C)$ in terms of $(2n+1)\times (2n+1)$ complex 
matrices which are skew-symmetric across the 
antidiagonal.  The triangular decomposition is the usual 
lower-diagonal-upper.  The involution $\Theta$ is inner, and 
given by 
$$\Theta =conj(\left(\matrix -1_{n\times n}&0&0\\
0&1&0\\
0&0&-1_{n\times n}\endmatrix \right)).$$
We have a diagram 
$$\matrix &&U=Spin(2n+1)\\
&\phi\nearrow&\\
U/K=S^{2n}&&pr\downarrow\\
&pr\circ\phi\searrow&\\
&&SO(2n+1)\endmatrix $$
The kernel of the projection is exactly the center of $U$, 
$\Bbb Z_2exp_{Spin}(i\pi h_n)$.  This central element is contained in 
$\phi (S^{2n})$, because $h_n$ corresponds to a $\Theta$-noncompact root, 
hence we can write 
$$exp_{Spin}(i\pi h_n)=i_{\alpha_n}(\left(\matrix 0&i\\
i&0\endmatrix \right))i_{\alpha_n}(\left(\matrix 0&i\\
i&0\endmatrix \right))^{-\Theta}.$$
(one can also observe that $\Theta$ must fix this element, 
since it is the unique nontrivial central element; note 
also that we can replace $\alpha_n$ by any of the $\Theta$-noncompact 
roots, and we get exactly the same central element).  
The central element acts (by multiplication) freely on 
the $\phi$-image, so that the image of $pr\circ\phi$ is isomorphic 
to $\Bbb R\Bbb P^{2n}$.  

Our calculations below are carried out inside of 
$SO(2n+1)$, and it is necessary to lift various maps to 
draw conclusions about $\phi (U/K)$.  This two-fold cover 
leads to some subtlties.  

We have 
$$i\frak p=\{X=\left(\matrix 0_{n\times n}&\matrix -\bar {z}_n\\
.\\
-\bar {z}_1\endmatrix &0_{n\times n}\\
z_n,..,z_1&0&\bar {z}_1,..,\bar {z}_n\\
0_{n\times n}&\matrix -z_1\\
.\\
-z_n\endmatrix &0_{n\times n}\endmatrix \right):z_j\in \Bbb C\}.$$
A Maple calculation implies that 
$$a(\frac {1-X}{1+X})=(1+2\sum_1^{n-1}\vert z_i\vert^2)^{h_1}..(1
+2\sum_1^1\vert z_i\vert^2)^{h_{n-1}}(1+2\vert z\vert^2)^{-(h_1+.
.+h_{n-1}+\frac 12h_n)}$$
$$=\prod_1^nv_j^{h_j},\quad v_1=\frac {1+2\sum_1^{n-1}}{1+2\sum_1^
n},..,v_n=\frac 1{(1+2\sum_1^n)^{1/2}}$$
$1>v_1\ge v_2\ge ..\ge v_{n-1}\ge v_n^2>0$.  Another miracle may occur 
in the calculation of volume.  

(here we did a calculation inside $SO(2n+1)$, and we have 
lifted the calculation).  All the entries here are positive.  
This implies that there is just one open component for 
projective space, the $pr\circ\phi$ image, but there are two 
open components for the sphere (they must be upper and 
lower hemispheres, so that the complement is $S^{2n-1}$).  

We must compute 
$$\int (1+\sum_1^{n-1})^{-i\alpha_1}(1+\sum_1^{n-2})^{-i\alpha_2}
..(1+\sum_1^1)^{-i\alpha_{n-1}}(1+\sum_1^n)^{-n-\frac 12+i(\alpha_
1+..+\alpha_n)}du_1..du_n,$$
where we set $u_i=2\vert z_i\vert^2$, so that $0\le u_i\,<\infty$ (Notes:  
$\alpha_1=\lambda_1-\lambda_2$$,..,\alpha_{n-1}=\lambda_{n-1}-\lambda_
n$$,\alpha_n=\lambda_n$; the $h_j$, $1\le j\le n-1$, 
are the usual, while $h_n$ has a $2$ in the $(n,n)$ position, and 
a $-2$ in the $(n+2,n+2)$ position; to calculate the volume 
element, we used the expression for the invariant metric 
on pages 13 and 14 of [Harmonic Analysis] - this is not 
an obvious calculation).  Integrating first with respect to 
$u_n$ (which only appears in the last factor), we see that 
$()$ equals 
$$\frac 1{(n-\frac 12-i\sum_1^n\alpha_j)}\int (1+\sum_1^{n-1})^{-
i\alpha_1}(1+\sum_1^{n-2})^{-i\alpha_2}..(1+\sum_1^1)^{-i\alpha_{
n-1}}(1+\sum_1^{n-1})^{-(n-1)-\frac 12+i(\alpha_1+..\alpha_n)}.$$
An induction now shows that this integral equals 
$$=\frac 1{(1-\frac 12-i\alpha_n)(2-\frac 12-i(\alpha_{n-1}+\alpha_
n))..(n-\frac 12-i(\alpha_1+..+\alpha_n))}$$
The $\Theta$-noncompact roots are the roots appearing in this 
formula (e.g.  $\lambda_1=\alpha_1+..+\alpha_n$ corresponds to the 
coordinate $z_n$,..,  $\lambda_1=\alpha_n$ corresponds to $z_1$).  Also 
$$\delta =(n-\frac 12)\lambda_1+(n-\frac 32)\lambda_2+..+(2-\frac 
12)\lambda_{n-1}+(1-\frac 12)\lambda_n$$
The formula above agrees with the conjectured formula.  

In the case of $S^{2n+1}$, to realize $\Theta$, we view $so(2n+2,
\Bbb C)$ 
as consisting of $3\times 3$ block matrices, where the block 
sizes are $n,2,n$.

$$i\frak p=\{X=\left(\matrix 0_{n\times n}&\matrix -\bar {z}_n&\bar {
z}_n\\
.&.\\
-\bar {z}_1&\bar {z}_1\endmatrix &0_{n\times n}\\
\matrix z_n,..,z_1\\
-z_n,..,-z_1\endmatrix &\matrix is&0\\
0&-is\endmatrix &\matrix -\bar {z}_1,..,-\bar {z}_n\\
\bar {z}_1,..,\bar {z}_n\endmatrix \\
0_{n\times n}&\matrix z_1&-z_1\\
.&.\\
z_n&-z_n\endmatrix &0_{n\times n}\endmatrix \right):z_j\in \Bbb C
\}.$$
A Maple calculation implies that 
$$a(\frac {1-X}{1+X})=\left(\matrix \endmatrix \right)$$

Now consider the case of $SU(2n)/Sp(n)$.  We can realize 
the involution in the following way:  given a $2n\times 2n$ 
matrix of complex numbers (of $trace=0$), we first 
reflect with respect to the antidiagonal, then we place 
minus signs against the two diagonal $n\times n$ blocks.  The 
fixed point set is our usual realization of $sp(n)$.  This 
involution is not inner.  

There are not any real roots (which for us is always 
the case because we are dealing with a maximal compact 
subalgebra $\frak h_0$ - see Knapp section on Cayley transforms).  
There are also just $n$ positive imaginary roots, 
corresponding to the root spaces along the antidiagonal; 
these imaginary roots are all of compact type.  Thus we 
are left to deal with the complex roots.  These come in 
pairs $\alpha ,\Theta (\alpha )$.  

In the cases in which $\Theta$ is outer, we need to supplement 
the noncompact roots with pairs $(\alpha ,\Theta (\alpha ))$ of complex 
roots.  
$$\frak p^{\Bbb C}=\bigoplus \frak g_{\alpha}\oplus\bigoplus_{(\alpha 
,\Theta (\alpha ))}\{X_{\alpha}-\Theta (X_{\alpha}):X_{\alpha}\in 
\frak g_{\alpha}\}$$
where the first sum is over the imaginary noncompact 
roots, and the second sum is over pairs of complex 
roots.  
  
Consider the case $Sp(n)/U(n)$.  In the case $n=2$ 
$$i\frak p=\{X=\left(\matrix 0&0&a_{13}&a_{14}\\
0&0&a_{23}&a_{13}\\
-\bar {a}_{13}&-\bar {a}_{23}&0&0\\
-\bar {a}_{14}&-\bar {a}_{13}&0&0\endmatrix \}\right)$$

We have 
$$a_{\phi}=(1+u_{32}-u_{41}-\vert detZ\vert^2)^{h_1}(1-2u_{31}-u_{
32}-u_{41}+\vert detZ\vert^2)^{h_2}$$
$$*(1+u_{32}-u_{41}-\vert detZ\vert^2)^{h_3}det(1+Z^{*}Z)^{-h_1-h_
2-h_3}$$

All $\bold w$ are possible:  
$$\bold w=1,\quad exp(i\pi h_2),\quad exp(i\pi h_1),\quad exp(i\pi 
(h_1+h_2)).$$
This appears to hold for all $n$.  

Suppose we consider $SO(2n)/U(n)$.  This case presumably 
does not have any noncompact type roots, so there 
should be a unique open component.

\bigskip

$==========================================$

\centerline{Appendix B. Poisson Geometry}

\bigskip

The rather uninspired calculations of $\S 1$ have a natural 
interpretation in terms of Poisson geometry.  The 
constructions involved are the same as in $\S 2$, and they 
are put into a beautiful general context in [EL].  We will 
be brief.  

To do calculations we will use the isomorphism of 
vector bundles 
$$U\times_Ki\frak p\to T(U/K):[u,x]\to\frac d{dt}\vert_{t=0}(ue^{
tx}K),$$
and we will use the Killing form to identify $\frak p^{*}$ with $
\frak p$.  

Consider the $Ad(T)$-stable Iwasawa decomposition of $\frak g$ as 
a direct sum of subalgebras:  
$$\frak g=\frak u\oplus (\frak n^{-}\oplus \frak h_0).\tag 2.15$$
Let $pr_{\frak u}$ denote the projection $\frak g\to \frak u$ along this 
decomposition.  Given $x\in \frak g$, with triangular 
decomposition $x=x_{-}+x_{\frak h}+x_{+}$, 
$$pr_{\frak u}(x)=(-x_{+}^{*}+x_{i\frak h_0}+x_{+}).\tag 2.16$$

The Evens-Lu Poisson bivector is given by 
$$\Pi ([u,x]\wedge [u,y])=\langle\Omega (u)(x),y\rangle ,\tag 2.17$$
where $\Omega (u):i\frak p\to i\frak p$ is given by 
$$\Omega (g_0)(x)=\{(pr_{\frak u}(ix^u))^{u^{-1}}\}_{i\frak p}.\tag 2.18$$
The operator $\Omega$ satisfies the equivariance condition 
$$\Omega (tuk)=Ad(k)^{-1}\Omega (u)Ad(k),\tag 2.19$$
for $t\in T$, $u\in U$, and $k\in K$.  

To understand $\Omega$, it is useful to consider the augmented 
operator $\tilde{\Omega }:\frak u\to \frak u$ given by 
$$\tilde{\Omega }(u)(x_{\frak k}+x_{i\frak p})=\{(pr_{\frak u}((x_{
\frak k}+ix_{i\frak p})^u))^{u^{-1}}\}.\tag 2.20$$
Relative to the decomposition $\frak u=\frak k\oplus i\frak p$, 
$$\tilde{\Omega }=\left(\matrix 1&*\\
0&\Omega\endmatrix \right).\tag 2.21$$

Now we need to factor $u=$ This is a huge mess.  

$===========================$

This augmented operator can be factored as the 
composition of four operators 
$$\frak g_0@>{I}>>\frak u@>{Ad(\bold u(g_0))}>>\frak u@>{T}>>\frak g_
0@>{Ad(\bold a_0^{-1}g_0)^{-1}}>>\frak g_0\tag 2.22$$
where the first operator is given by $I(x_{\frak k}+x_{\frak p})=
x_{\frak k}+ix_{\frak p}$, 
and $T(g_0)$ (which is reminiscent of a Toeplitz operator) 
maps $x=-x^{*}_{+}+(x_{t_0}+x_{i\frak a_0})+x_{+}$ to 
$$T(g_0)(x)=pr_{\frak g_0}(x^{\bold l'\bold a_1(g_0)})=$$
$$[(x_{+}^{\bold l'\bold a_1})_{+}]^{\sigma}+(x_{t_0}+(x^{\bold l'
\bold a_1}_{+})_{\frak t_0+\frak a_0})+(x_{+}^{\bold l'\bold a_1}
)_{+},\tag 2.23$$
where $\bold l'=\bold a_0\bold l\bold a_0^{-1}$, and the last equality depends upon the 
fact that conjugation by $\bold l'\bold a_1(g_0)$ maps $\frak n^{
-}$ into itself, and 
that $\frak n^{-}$ terms disappear when we use $(2.16)$.  
 
\proclaim{ (2.24) Lemma} (a) $\Omega\in so(\frak p)$ 

(b) $ker(\Omega (g_0))=\{[g_0,(\frak a_0^{\bold u(g_0)^{-1}})_{\frak p}
]\}$ 

(c) $Pfaffian(\Omega (g_0)\vert_{ker(\Omega )^{\perp}})=\bold a_1
(g_0)^{2\delta}$.  

(d) 
$$d\Omega\vert_{g_0}(x)=P_{\frak p}\circ [pr_{\frak g_0}^{Ad(g_0)^{
-1}},ad_x]$$

\endproclaim

\demo{Proof of (2.24)} For (a) let $X=x^{g_0}$, $Y=y^{g_0}$, 
$x,y\in \frak p$.  Then 
$$\kappa (\Omega (g_0)x,y)=\kappa (pr_{\frak g_0}(iX),Y)$$
$$=\kappa (-iX_{+}^{\sigma}+iX_{+},Y_{+}^{\sigma}+Y_{\frak h_0}+Y_{
+})=2\kappa (iX_{+},Y_{+}^{\sigma}).\tag 2.25$$
This is clearly skew-symmetric in $X$ and $Y$, because $\sigma$ 
preserves $\kappa$ and it is complex antilinear.  

For (b), note that $(2.23)$ implies the kernel of $T$ is $i\frak a_
0$.  
Thus $(2.22)$ implies 
$$ker(\tilde{\Omega }(g_0))=\{[g_0,x]:(x_{\frak k}+ix_{\frak p})\in 
i\frak a_0^{\bold u(g_0)^{-1}}\}\tag 2.26$$
This, together with $(2.21)$, implies $(b)$.  

For (c), note that in $(2.22)$ the first, second and fourth 
operators preserve volume determined by the Killing 
form.  The determinant of $T$ (relative to the Killing 
form volumes) is the same as the determinant of the 
operator on $\frak n^{+}$ which maps $x_{+}$ to $(x_{+}^{\bold l'
\bold a_1})_{+}$.  This 
determinant equals 
$$\prod_{\alpha >0}\bold a_1^{2\alpha}=\bold a_1^{4\delta}\tag 2.27$$
Thus the Pfaffian is $\bold a_1^{2\delta}$.  

(d) is straightforward.  \qed
\enddemo

By $(b)$ the tangent directions in $G_0/K$ determining the 
symplectic leaves are given by $[g_0,x]$ such that $x^{\bold u}\perp 
\frak a_0$.  
This is clearly $A_0$-invariant, because $\bold u(a_0g_0)=\bold u
(g_0)$.  
Thus the left action of $A_0$ permutes the symplectic 
leaves.  The symplectic form is given by the formula 
$$\omega ([g_0,x],[g_0,y])=\kappa (\Omega (g_0)\vert_{ker(\Omega 
)^{\perp}})^{-1}(x),y)\tag 2.28$$
Since $A_0$ permutes the symplectic leaves, we can regard 
this as defining (in an indirect way) a symplectic form 
on the quotient $A_0\backslash G_0/K$ [There has got to a more 
elegant way of thinking about this].  

\proclaim{ (2.29) Proposition}(a) The action of $T_0$ is 
Hamiltonian with momentum map 
$$\mu :A_0\backslash G_0/K\to (\frak t_0)^{*}:A_0g_0K\to\langle i
log(\bold a_1(g_0)),\cdot\rangle ,$$
This momentum map is proper, and it is semibounded.  

(b) The symplectic measure is 
$$\omega^d/d!=\bold a_1(g_0K)^{-2\delta}dV_{A_0\backslash G_0/K}(
A_0g_0K)$$
(where the invariant measure is suitably normalized).  

\endproclaim

Part (a) is proven in [FO].  Part (b) follows from (c) of 
$(2.24)$.  

We can now apply the Duistermaat-Heckman exact 
stationary phase, as generalized to noncompact manifolds 
as in [PW].  

Suppose that $A_0g_0K$ is fixed by $T_0$.  Then $g_0$ normalizes 
$T_0$, and this implies $a_0g_0\in K$ for some $a_0\in K$.  Thus 
there is just one $T_0$ fixed point, the basepoint.  If $X$ 
denotes the element of $\frak t_0$ corresponding to $\delta +\Lambda$, then 
the Pfaffian of the infinitesimal action of $X$ at the 
basepoint equals 
$$Pf(ad(X)\vert_{\frak p})=\prod\langle\delta +\Lambda ,\alpha\rangle 
,\tag 2.30$$
where the product is over pairs of positive roots 
$(\alpha ,\Theta (\alpha ))$ which are not of compact type.  This concludes 
the proof of $(2.2)$.  \qed

\bigskip

\centerline{Comments and Questions}

\bigskip

1.  Consider the diagram $(0.3)$.  This diagram exists for 
any automorphism $\Theta$ of $U$.  However, the maps are 
totally geodesic, and there is a simple algebraic 
characterization of the images, only when $\Theta$ is 
involutive.

2.  As in [EL], let $\Cal L(\frak g)$ denote the Grassmannian 
consisting of (real) subalgebras of $\frak g$ which are Lagrangian 
with respect to the imaginary part of the Killing form 
(for example any real form of $\frak g$ is a point in $\Cal L(\frak g
)$).  The 
group $G$ acts on $\Cal L(\frak g)$ through the adjoint action.  

Modulo covering space issues, the diagram $(0.3)$ has an 
interpretation in this context:  
$$\matrix Ad(U)\cdot \frak g_0&@>{\phi}>>&Ad(U)\\
\downarrow&&\downarrow\\
Ad(G)\cdot \frak g_0&@>{\phi}>>&Ad(G)&{{\psi}\atop {\leftarrow}}&
Ad(G)\cdot \frak u\\
&&\uparrow&&\uparrow\\
&&Ad(G_0)&{{\psi}\atop {\leftarrow}}&Ad(G_0)\cdot \frak u\endmatrix 
.\tag $0.3'$ $$
The interesting question is whether this picture can be 
given a coherent Poisson-theoretic interpretation.  

On the one hand if we consider the bialgebra structure 
$(\frak u,\frak n^{-}\oplus \frak h_0)$ for $\frak g$, then we can equip $
G$ and $U$ with Lie 
Poisson structures, and  $Ad(U)\cdot \frak g_0$ and $Ad(G)\cdot \frak g_
0$ inherit 
homogeneous Poisson structures.  On the other hand if 
we consider the bialgebra structure $(\frak g_0,\frak n^{-}\oplus 
i\frak h_0)$ for $\frak g$, 
then we can equip $G$ and $G_0$ with Poisson structures, 
and $Ad(G)\cdot \frak u$ and $Ad(G_0)\cdot \frak u$ inherit homogeneous Poisson 
structures.  

There is an obvious consistency question:  is the Poisson 
structure on $Ad(G)$ the same in both cases?  One 
suspects that the answer is no, but I am uncertain about 
this.  

We can now ask a number of questions (the precise 
meaning of which may be dependent upon the resolution 
of the consistency question above).  

(a) Are the maps in $(0.3')$ Poisson maps?  

(b) Assuming the truth of (a), does this shed some light 
on the intersection of the embeddings with the 
triangular decomposition for $G$?  

(c) What has this got to do with the Poisson structure 
on $U/K$ which was introducted by Foth and Lu?  Is this 
the same as the Poisson structure introduced above?  
Tentatively I think the answer is no.  Their structure 
apparently depends upon the choice of a triangular 
decomposition which intersects nicely with $\frak g_0$.  And yet 
they seem to be obtaining the same sort of 
stratification.  What is the momentum map in their 
setup?  Maybe there are multiple ways to obtain 
momentum maps so that the Duistermaat-Heckman result 
can be applied.  This seems too bizarre.  

(d) There must be explanation for the equality of 
integrals in terms of this picture $(0.3')$.  For example 
possibly there is some kind of deformation which allows 
these integrals to swim around in $G$.  

3.  If we do not assume that $\Theta$ is involutive, can we 
still say something about the integrals?  

4.  We can also consider replacing the complex group $G$ 
by a real group $G$.  The structures can be modified 
appropriately to obtain a similar picture.  Is there a 
Poisson interpretation?  Can we say something about 
integrals?  Can this be related to harmonic analysis?  
For example to more general $\bold c$-functions.

5.  In the text we used a Poisson structure on $G_0/K$.  Is 
this the same as the Poisson structure coming from the 
Foth-Lu construction, or are these different?  

6.  Can we write down this Poisson structure in an 
explicit way?  

7.  Is there a geometric explanation, based upon general 
principles, for the fact that the $\phi$-images intersect the 
triangular decomposition in a way which preserves 
homotopy equivalence?  

\bigskip

$==========================================$ 

$=================$ 

Many examples are worked out in the Cayley coordinates 
section.  At this point the pattern is not clear.  Some 
examples:  In the case of $Gr(n,\Bbb C^{n+m})$, the number of $-1$'s 
in the upper $n\times n$ block must equal the number of $-1$'s in 
the lower $m\times m$ block; this condition intuitively arises 
from the fact that for each noncompact root, $exp(i\pi h_{\alpha}
)$ 
has a $-1$ in each of these blocks.  For a sphere (at least 
of even dimension, inner), there are two open 
components; this case is unusual, in that $\bold w=exp(i\pi h_{\alpha}
)$ 
is the same, the central element in the spin group, for 
all noncompact type roots.  The same may be true for 
the odd dimensional, outer case.  For $U(2n)/Sp(n)$ 
(outer), there are not any noncompact roots, but there 
are complex roots, so there are multiple components.  
For $Sp(n)/U(n)$ (inner), anything in $T_0^{(2)}$ is possible (this 
seems to be the most straightforward example; it's 
interesting that this includes my candidate for string 
theory space).  For $SO(2n)/U(n)$, the number of $-1$'s in 
the upper left corner must be even (and equal to the 
number of $-1$'s in the bottom right corner).  It is 
surprising that even these Hermitian symmetric cases 
seem somewhat eclectic.  Other Hermitian symmetric 
cases:  $SO(p+2)/S(O(p)\times O(2))$, $E_6/(SO(10)\times \Bbb T)$, 
$E_7/(SO(12)\times \Bbb T)$.  Other classical cases:  $SU(n)/SO(n
)$ 
(outer); one might need to consider even $n$ and odd $n$ 
separately.  $SO(n+m)/S(O(n)\times O(m))$ (outer, if $n$ or $m$ is 
odd) and $Sp(n+m)/Sp(n)\times Sp(m)$.  I do not have much 
insight into the exceptional cases (listed on page 518 of 
[Helgason]).  

$=================$

It is natural to ask whether there is a more general 
formula for an integral (corresponding to lower strata) 
of the form 
$$\int_{\Sigma_{\bold w}^{\phi (U/K)}}a_{\phi}(\bold w_1g\bold w_
1^{-\Theta})^{-i\lambda}$$
where $\bold w=\bold w_1\bold w^{*\Theta}_1$, with respect to the natural 
background Riemannian measure.  

More generally one could ask whether there is a 
generalization involving automorphisms $\Theta$ which are not 
necessarily involutive.

\bigskip

\centerline{Comments on Weyl dimension formula approach}

\bigskip

The Harish-Chandra formula in the introduction can be 
deduced from the Weyl dimension formula (see the 
introduction to [Pi2]).  This naturally raises the question 
whether the same approach applies in this more general 
context.  

Suppose that we consider the special value $-i\lambda =\Lambda$, 
where $\Lambda /2$ is dominant integral.  For simplicity we 
suppose that we are in the inner case).  We can write 
$$\int_{\phi (U/K)}a_{\phi}^{-i\lambda}=\int_U\vert\sigma_{\Lambda 
/2}(gg^{-\Theta})\vert^2dg=\int_U\vert\langle g^{\Theta}vac,gvac\rangle
\vert^2dg$$
where $\sigma_{\Lambda /2}$ and $vac$ denote the fundamental matrix 
coefficient and normalized vacuum correspond to the 
highest weight module $V=V_{\Lambda /2}$ corresponding to $\Lambda 
/2$.  
There is $U$-equivariant linear mapping 
$$V\otimes\bar {V}\to L^2(K\backslash U):v\otimes\bar {w}\to\langle 
g^{\Theta}v,g\bar {w}\rangle ,$$
[or, in terms of operators, $L\to trace(g^{\Theta}Lg^{-1})$], where $
U$ 
acts on the tensor product by $g:v\otimes\bar {w}\to g^{\Theta}v\otimes 
g\bar {w}$ [or, in 
terms of operators, $g:L\to g^{\Theta}Lg^{-1}$].  We can decompose 
$$L^2(K\backslash U)=\bigoplus W$$
as a multiplicity free representation, where the sum is 
over irreducible representations of $U$ which contain a 
(essentially unique) $K$-invariant vector.  Therefore the 
map $()$ essentially projects $V\otimes\bar {V}$ to its spherical 
subrepresentation, which will necessarily be multiplicity 
free.  

The projection to the $K$-invariant vectors is given by 
averaging over $K$.  At some point the $\bold w$ should come 
into play.  One way:  if we consider the vacuum, the 
corresponding projection $P_{vac}$ is fixed by the $\bold w$.  By 
why does this give rise to the sum over $\bold w$.  This is 
baffling.

\demo{Proof} We will consider a slight reformulation of 
the problem.  We identify $U/K$ with $\phi (U/K)$.  Let $S$ 
denote the inverse image of $\Sigma_1^{U/K}$ in $U$, with respect to 
the projection $U\to U/K$.  Consider the diffeomorphism 
$$\Psi :exp(i\frak a_0)\times A_0\backslash G_0\to S:(t,A_0g_0)\to 
t^{-1}\bold u(g_0).$$
We must show that the Jacobian for the mapping $\Psi$, 
with respect to the Riemannian structures induced by 
the Killing form, is equal to a constant times $a_{\phi}^{2\delta}$.  To 
do this we identify $i\frak a_0$, $\frak g_0\ominus \frak a_0$, and $
\frak u$ with the tangent 
spaces to $exp(i\frak a_0)$, $A_0\backslash G_0$ and $U$, respectively, using the 
exponential map and an appropriate translation.  Let 
$P:\frak g\to \frak u$ denote the projection with kernel $\frak n^{
-}\oplus \frak h_{\Bbb R}$.  We 
compute 
$$d\Psi\vert_{(t,A_0g_0)}:i\frak a_0\oplus (\frak g_0\ominus \frak a_
0)\to \frak u:(\chi ,x)\to\frac d{d\epsilon}\vert_{\epsilon =0}(t
e^{\epsilon\chi})^{-1}\bold u(g_0e^{\epsilon x})\bold u(g_0)^{-1}
t$$
$$=Ad(t^{-1})\{-\chi +P(Ad(\bold u)(x))\}.$$

The operator $Ad(t^{-1})$ preserves $\frak u$-volume, so it can be 
ignored.  

Now suppose that we first consider the inner case, i.e.  
$\frak a_0=0$.  Because $\bold u(g_0)=\bold a^{-1}\bold l^{-1}g_0$, the derivative can be 
expressed as a composite map 
$$\frak g_0@>{Ad(g_0)}>>\frak g_0\to \frak u:x\to Ad(g_0)(x)\to P
(Ad(\bold a^{-1}\bold l^{-1}g_0)(x))$$
Because $G_0$ is semisimple, $Ad(g_0)$ preserves volume.  In 
the general case we still need to worry about $g_0$, 
because we have the restricted domain $\frak g_0\ominus \frak a_0$.  

Write $\bold a=\bold a_1\bold a_0$, relative to the decomposition 
$A=exp(i\frak t_0)A_0$.  Since $\bold a_0\in G_0$, $Ad(\bold a_0)$ will also preserve 
volume.  Thus the determinant of $()$ equals the 
determinant of the map 
$$\frak g_0\ominus \frak a_0\to \frak u\ominus i\frak a_0:x\to P(
Ad(\bold l')Ad(\bold a_1^{-1})(x)),$$
where $\bold l'=\bold a_1^{-1}\bold la_1\in N^{-}$.  

Given $x\in \frak g$, we write $x=x_{-}+h+x_{+}$ for its triangular 
decomposition.  The $\frak n_{-}\oplus \frak h_{\Bbb R}\oplus \frak u$ decomposition is given by 
$$(x_{-}+x_{+}^{*})\oplus\frac 12(h+h^{*})\oplus (-x_{+}^{*}+\frac 
12(h-h^{*})+x_{+}).$$
If $x\in \frak g_0$, then $x_{-}=x_{+}^{\sigma}$, and $h=t_0\in \frak t_
0$.  Because $\bold l'\bold a^{-1}$ 
maps $\frak n_{-}$ into itself, we see that $()$ equals 
$$P(Ad(\bold l'\bold a_1^{-1})(x_{+}^{\sigma}+t_0+x_{+}))=-[(x_{+}^{
\bold l'\bold a_1^{-1}})_{+}]^{*}+t_0+(x_{+}^{\bold l'\bold a_1^{
-1}})_{+}.$$
Thus the determinant of $()$ is the same as the (real) 
determinant of the map $x_{+}\to (x_{+}^{\bold l'\bold a_1^{-1}})_{
+}$.  Because of the 
unipotence of $Ad(\bold l')$, this is equal to 
$$\prod_{\alpha >0}\vert \bold a_1^{-\alpha}\vert^2=\bold a_1^{-4
\delta}=a_{\phi}^{2\delta}$$

This completes the proof in the inner case.  It is 
possible that the general case involves some refinement 
which our formula does not account for.  

\qed
\enddemo

(note that in our semisimple context, the $+1$ eigenspace 
is the orthogonal complement of the $-1$ eigenspace, 
relative to the Killing form, so that $Ad(g)\circ\sigma$ is 
determined in a very explicit way by its $-1$ eigenspace).